\documentclass[11pt]{article}

\usepackage[latin1]{inputenc}
\usepackage[T1]{fontenc}

\usepackage{rotating}
\usepackage{multirow}
\usepackage{slashbox}
\usepackage{srcltx}
\usepackage{amsmath,latexsym}
\usepackage{pstricks}
\usepackage{pst-node}
\usepackage{amssymb}
\usepackage{epsfig}
\usepackage{fancyhdr}
\usepackage{subfigure}
\usepackage{xspace}
\usepackage{enumerate}
\usepackage{rotating}
\usepackage{multirow}
\usepackage{pst-node,colortbl}

\newtheorem{ejem}{\bf Example}[section]
\newtheorem{rem}{\bf Remark}[section]

\newcommand{\tl}{\tilde \lambda}

\allowdisplaybreaks[1]

\begin{document}
\title{Ordered Median Hub Location Problems with Capacity Constraints}

\author{J. Puerto$^a$ \and A.B. Ramos$^a$ \and A.M. Rodr{\'\i}guez-Ch{\'\i}a$^b$ \and
M.C. S\'anchez-Gil$^b$ }

\date{
\small
$^a$ Facultad de Matem{\'\i}ticas. Universidad de Sevilla, Spain.\\
$^b$ Facultad de Ciencias. Universidad de C\'adiz, Spain.\\
\today
}

\maketitle \thispagestyle{empty}
\begin{abstract}
The Single Allocation Ordered Median Hub Location problem is a recent hub model introduced
in \cite{PRRCH11} that provides a unifying analysis of  a wide class of hub  location models.
In this paper, we deal with the capacitated version of this problem, presenting two formulations as well as some preprocessing phases for fixing variables. In addition, a
strengthening of one of these formulations is also studied through the use of some families of valid inequalities.
A battery of test problems with data taken from the AP library are solved where it is shown that
the running times have been significantly  reduced with the improvements presented in the paper.

\end{abstract}

\section{Introduction}\label{Sec:Intro}

Network design problems are among the most interesting models in combinatorial optimization. In the last years researchers have devoted a lot of attention to a particular  member within this family, namely the hub location problem, that combines network design and location aspects of supply chain  models,   see the surveys \cite{SB08,CEK02,COK12}.
 The main advantage  of using  hubs in distribution problems is that they allow to consolidate shipments in order to reduce transportation costs by applying  economies of scale; which are naturally incorporated to the models  through discount factors.
Hub location problems have been studied from different perspectives   giving rise to a number of papers  considering different criteria to be optimized: the minimization of the overall transportation cost (sum)  (see
\cite{CGM07,C96,EK99,LYG05,M05,M05b,MCL06}), the minimization of the largest
transportation cost or the coverage cost (\cite{BLR06,CLZ07,KT00,KT03,KS06,MEK09,TK07,W08}), et cetera.

Apart from the choice of the optimization criterion,  another crucial aspect in the literature {on} hub location, and in general {on} any location problem, is the assumption of capacity constraints. One can recognize that  although
this assumption implies  more realistic models, the difficulty to solve them also increases in
orders of magnitude with respect to {their} uncapacitated counterpart.
In many cases new formulations are needed and a more specialized analysis is often required to solve
even smaller sizes than those previously addressed for the uncapacitated versions of the problems.
For this reason, capacitated versions of hub location problems have attracted the interest of locators
in the last years, see  \cite{A94,CAMP94b,CDF09,CNS10,CNS10b,EKEB00,BKEE04,EK99,M05}.  In the same line, we also mention some other references related with congestion at hubs, as congestion acts as a limit  on capacity, see \cite{CMFL09,EH04,MS03}.

An interesting version of  hub location model is the Capacitated Hub Location Problem with
Single Allocation  (CSA-HLP), see {\cite{CDF09,CNS10,EK99}}.
In this  context, single allocation means that  incoming and outgoing flow of each site must be shipped via the same hub. In contrast to single allocation models, where binary variables are required in the allocation phase,  multiple allocation allows different delivery patterns which in turns implies the use of continuous variables simplifying  the problems.
The CSA-HLP model incorporates capacity constraints on the incoming flow at the
hubs coming  from origin sites 
or even simpler, on the number of non-hub nodes assigned to
each hub.
The inclusion of capacity constraints make these models challenging from a theoretical
point of view. Regarding its applicability   we cite  one example described in
Ernst at al. \cite{EK99} based on a postal delivery application, where a set of $n$ postal districts
(corresponding to postcode districts  represented by nodes) exchange daily mail. The mail
between all the pairs of nodes must be routed via one or at most two mail consolidation centers
(hubs). In order to meet time constraints, only a limited amount of mail could be sorted at each
sorting center {(mail is just sorted once, when it arrives to the first hub from origin sites)}. Hence, there are capacity restrictions on  the incoming mail that must be sorted. The
problem requires to choose the number and location of hubs, as well as  to determine the
distribution pattern of the mail.

The CSA-HLP has received less attention in the
literature than its uncapacitated counterpart. Campbell \cite{CAMP94b} presented the first
integer  Mathematical Programming formulation for the Capacitated Hub Location Problem.
This formulation was strenghthened by  Skorin-Kapov et al. \cite{SKO96}.
Ernst and Krishnamoorthy \cite{EK99}, 
 proposed a new model involving {three-index}
continuous variables  and  developed a solution approach based on Simulated Annealing 
where the bounds obtained are embedded in a branch-and-bound procedure devised for solving the
problem optimally. Recently, Correia et al.  \cite{CNS10} have shown that this formulation may
be incomplete and an additional set of inequalities is proposed to assure the validity of the model in all situations. A new formulation using only two indices variables was proposed by
Labb\'e et al. \cite{LYG05}, where a polyhedral analysis and  new valid inequalities were
addressed. Although  this formulation has only a quadratic number of variables, it has an exponential number  of  constraints, and to solve it the authors  developed a branch-and-cut algorithm based on {their} polyhedral analysis. Contreras et al.  \cite{CDF09} presented for the same problem  a Lagrangean relaxation enhanced with reduction tests that allows the computation of tight upper and lower bounds for a large set of instances.

In two recent papers, \cite{PRRCH11,PRRCH13}, a new model of hub location, namely the Single
Allocation Ordered Median hub location problem (SA-OMHLP), has been introduced and analyzed.
This problem can be seen as a powerful tool from a modeling point of view since it allows a common framework to represent many of the previously considered criteria in the literature of hub location. Moreover, this approach is a natural way to represent the
differentiation of the roles played by the different parties (origins, hubs and destinations) in logistics networks  \cite{FGOS10,KNPR10,KNPR09,MNPV08}.
This model does not assume, in advance,  any particular structure on the network (\cite {CDF09,CDL11}). Instead of that this structure is derived from the choice of the parameters defining the objective function.
Apart from the above mentioned characteristics, ordered median objectives are also useful to obtain robust solutions in hub problems by applying $k$-centrum, trimmed-mean or anti-trimmed-mean criteria.
It is worth mentioning that although it is called single allocation,  its meaning  slightly differs from the classical interpretation in hub location  where  each site is allocated to just one hub and all the  incoming and outgoing flow to-from this site is shipped  via the same link (the one joining this site and its allocated hub). In this model, single allocation means that all the outgoing flow is delivered through the same hub, but the incoming flow can come from different hubs. Actually, this is a mixed model and basically the same situation described above, about postal deliveries, naturally fits in this framework assuming that  letters from the same origin should be sorted, with respect to their destinations, in the same place and from there they are delivered via their cheapest routes. Observe that in this scheme it is { also} natural that incoming flow in a final destination comes from different hubs.

The SA-OMHLP  distinguishes among segmented origin-destination deliveries giving different scaling factors to the origin-hub, hub-hub and hub-destination links.
The cost of each origin-first hub link is scaled by a factor that depends on the position of this cost in the ordered sequence of costs from each origin to its corresponding first hub \cite{BDNP05,MNPV08,NaP05}. Moreover, the overall interhub cost and hub-destination cost are multiplied by other economy of scale factors. The goal is to  minimize the overall shipping cost under the above weighting scheme. The reader may note that
the first type of scaling factors mentioned above adds a ``sorting'' problem to the underlying  hub location model,  making its formulation and solution much more challenging.
This model and two different formulations were introduced in  \cite{NaP05}  while a specialized B\&B\&Cut algorithm was developed in \cite{PRRCH11,PhDAna}.
None of those formulations could handle capacities since the computation burden of the problems were highly demanding.   Thus, the SA-OMHLP with capacity constraints, {i.e. Capacitated Single Allocation  Ordered Median hub location problem (CSA-OMHLP)} is currently an open line of research for further analysis.

In this paper, we analyze in depth the CSA-OMHLP trying to obtain a better knowledge and alternative ways to solve it.
Thus, the contributions of this paper are threefold. First, it combines for the first time three challenging elements in location analysis: hub facilities, capacities and ordered median objectives; proposing a promising IP formulation which remarkably reduces the number of decision variables.
Second, this paper  strengthens that formulation with variable fixing and some families of valid inequalities that have not been considered before. Finally, despite the difficulty of considering simultaneously capacitated models, hubs and ordering, the techniques proposed in this paper allow to solve instances of similar sizes to those already considered in the literature for simpler models (uncapacitated and multiple allocation \cite{KNPR10}).

The paper is organized as follows: {in Section \ref{Sec:Model}} we will provide,  first, a MIP formulation for the capacitated version of the problem extending the one in \cite{PRRCH11}  and then another formulation  in the
spirit of \cite{PRRCH13}  where the number of variables has  been considerably reduced with respect
to the previous one.
Section \ref{Sec:StrenForm} strengthens the latter formulations with variable fixing
and  several new families of  valid inequalities. 
{In Section \ref{Sec:Results}}, the effectiveness of the proposed methodology is shown with an {extensive} computational experience comparing the performance of the two formulations and the strengthening proposed along the paper.
Finally, the paper ends with some conclusions.

\section{Model and MIP formulations}\label{Sec:Model}

The goal of this paper is to analyze the CSA-OMHLP. For this reason, we elaborate from the most promising formulations of the non-capacitated version of that problem, namely the so called radius
(covering) formulations, see   \cite{PRRCH11, PRRCH13}. In order to be self-contained and
for the sake of  readability, we   include  next a concise description
of these formulations in their application to the capacitated problem.

Let  $A=\{1,\ldots,n\}$ be a set of $n$ client sites, where each site is collecting or gathering some commodity that must be  sent to the remaining ones.  It is assumed, without loss of
generality, that the set of candidate sites for establishing hubs is  also A.
Let $w_{jm} \ge 0$ be the amount of commodity to be  supplied from the  $j$-th to the $m$-th site for all $j,m \in A$, and  let $W_j=\sum_{m \in A} w_{jm}$ be the total amount of commodity to be sent from  the $j$-th
site. Let $c_{jm}\ge  0$ denote the unit cost of sending commodity from site $j$ to site $m$ (not
necessarily satisfying the triangular inequality). It is assumed free self-service, { i.e.,} $c_{jj}=0$, $\forall j\in A$.

Let $p \le  n$ be the number of hubs to be located and 
let $b_j$ be the capacity of a hub  located at site $j$, with $j \in A$.
A solution for the problem is a set of sites $X\subseteq A$ with $|X| = p$ and enough capacity to cover the flow coming from the sites; plus a set of  links connecting pairs (flow patterns) of sites $j$, $m$ for all $j,m \in A$. Moreover, it is assumed that   the flow pattern between each pair of sites traverses at least one and no more than  two
hubs from $X$.

As it was mentioned in Section \ref{Sec:Intro}, {the main advantage of using hubs is to reduce costs by applying economies of scale to consolidated flows in some part of the network. In this model} the transportation cost is decoupled into the three differentiated possible links: origin site-first hub, hub-to-hub, and hubs-final destination. These transportation costs are  scaled in a different way.
The model weights origin site-first hub transportation costs  by using
parameters $\lambda=(\lambda_1,\ldots,\lambda_n)$, with $\lambda_i \ge 0 \,\, \forall i \in A$, depending on their {ordered rank values}.
This is, let  $\hat c_{jk}$ be the cost of the overall flow sent from the origin site $j$ if it were delivered via the first hub $k$, i.e.
$ \hat c_{jk}:= c_{jk}W_j, \, j,k\in A.$ Next, if $\hat c_{jk}$ were ranked in the $i$-th position among all these costs,  then  this term would be scaled by $\lambda_i$ in the objective function.
For the two remaining links there is a compensation factor $0<\mu<1$ for the deliveries between hubs, and another one $0<\delta<1$, $\mu < \delta$, for the deliveries between hubs and final
destination sites. These parameters  may imply that, even in the case where the costs satisfy the triangle inequality,  using a second  hub results in a cheaper connection than going directly from the first hub to the final destination.  Actually, it represents the application of the economy of scale by the consolidation of flow in the hubs.


In the following we present a first valid formulation of the 
CSA-OMHLP, based on covering variables (the reader is referred to
\cite{PRRCH11, PRRCH13} for further details).
Sorting the different delivery costs values  $(\hat c_{jk})$ for $j,k\in A$, in increasing order, we get the ordered cost sequence:
$$\hat  c_{(1)}:= 0 <\hat  c_{(2)} < \cdots <\hat
c_{(G)} := \max\limits_{1 \le j,k \le n} \{\hat{c}_{jk}\}.$$
where $G$ is the number of different elements of the above cost sequence. For convenience we consider $\hat{c}_{(0)}:=0$.

For $i \in A$ and $h = 1,\ldots,G$, we define the following  set of covering variables,
\begin{equation} \label{xjkdef}
\bar{u}_{ih} := \left\{ \begin{array}{cl}
                  1, & \mbox{if the $i$-th smallest allocation cost is at least
                  $\hat{c}_{(h)}$,}\\[2ex]
                  0, & \mbox{otherwise.}
                  \end{array} \right.
\end{equation}
Clearly, the $i$-th smallest allocation cost is equal to
$\hat{c}_{(h)}$ if and only if $\bar{u}_{ih} = 1$ and $\bar{u}_{i,h+1} = 0$. \\
 In addition, this formulation uses two more sets of variables:
\begin{eqnarray}
\label{var:r}
x_{jk} &=& \begin{cases}
1, & \mbox{if the commodity sent from origin site $j$ goes first to the hub $k$,}\\
0, & \mbox{otherwise.}
\end{cases} \\
\nonumber
 s_{k\ell m} &=& \mbox{flow that goes  through  a first hub $k$ and a second  hub $\ell$ with destination
 $m$},
\end{eqnarray}
with $j,k,\ell,m\in A.$ Since we assume that any origin is allocated to itself if it is a hub,
the above definition implies that site $k$ is opened as a hub if the corresponding variable $x_{kk}$ {takes the value 1}.

The  formulation of the model is:
\begin{eqnarray}
\label{CCOV0}
{(F_{\bar{u}})}&\min & \sum_{i \in A} \sum_{h=2}^{G} \lambda_{i}(\hat c_{(h)}-\hat c_{(h-1)})\bar{u}_{ih}+
\sum_{k \in A} \sum_{\ell \in A} \sum_{m \in A} (\mu c_{k \ell}+\delta c_{\ell m})s_{k\ell m}\\
&s.t. &
\label{CCOV1}
\sum_{k \in A}
x_{jk}=1, \quad \forall j\in A \\
\label{CCOV2}
& &\sum_{j \in A}  x_{jk} \leq nx_{kk}, \quad \forall k\in A\\
\label{CCOV10}
& &\sum_{k \in A} x_{kk}=p\\
\label{CCOV5}
& &\sum_{\ell \in A} s_{k\ell m}=\sum_{j \in A} w_{jm}x_{jk} , \quad \forall k,m\in A \\
\label{CCOV4}
& & s_{k\ell m}\leq  \sum_{j \in A} w_{jm} (1-x_{mm}) \quad \forall k,\ell,m\in A, \quad \ell\neq m\\
\label{CCOV6}
& &\sum_{\ell \in A}\sum_{m \in A} s_{k\ell m}\leq x_{kk}\sum_{j \in A} W_{j} , \quad \forall k\in A \\
\label{CCOV7}
& &\sum_{k \in A}\sum_{m \in A} s_{k\ell m}\leq x_{\ell \ell}\sum_{j \in A} W_{j} , \quad \forall \ell\in A \\
\label{UUVC9}
& & {   \sum_{j \in A} W_ j x_{jk}\leq b_kx_{kk}, \quad \forall k\in A} \\
\label{CCOV8}
& & \sum_{i \in A} \bar{u}_{ih}= \sum_{j \in A}
\sum_{\stackrel{k=1}{\hat{c}_{jk}  \ge \hat  c_{(h)}}}^n  x_{jk}, \quad \forall h=1,\ldots,G \\
\label{CCOV9}
& & \bar{u}_{ih}\ge \bar{u}_{i-1,h}, \quad \forall i\in A\setminus\{1\}, \: h=1,\ldots,G\\
\label{CCOV11}
& & \bar{u}_{ih},x_{jk} \in \{0,1\},  s_{k \ell m} \ge 0, \quad \forall i,j,k, \ell,m\in A, \quad h=1,\ldots, G
\end{eqnarray}

The objective function \eqref{CCOV0} accounts for the weighted sum of the three components of the shipping cost, namely origin-first hub, hub-hub and hub-destination.
The origin-hub costs are  accounted after  assigning the lambda parameters, i.e.
 $\sum\limits_{i \in A}
\sum\limits_{h=2}^G \lambda_{i} \cdot (\hat c_{(h)} -\hat  c_{(h-1)})
\cdot \bar{u}_{i h}.$
In addition,  the second and third blocks of delivery costs, i.e. the hub-hub and hub-destination cost,
 scaled with the $\mu$ and $\delta$
parameters respectively, can be stated as:
$\sum_{k \in A} \sum_{\ell \in A} \sum_{m \in A} (\mu c_{k \ell}+\delta c_{\ell m})s_{k\ell m}.
$

The constraints of the model are described in \cite{PRRCH11} with the only exception of (\ref{UUVC9}), that is the capacity constraint on the hubs, 
{in spite of that and for the sake of completeness  we include below a brief description of them.}
{Constraints  (\ref{CCOV1}) ensure that  the flow from  the origin site $j$ is associated with a unique first hub.
Constraints (\ref{CCOV2}) ensure that any origin only can be allocated to an open hub.
Constraint (\ref{CCOV10}) fixes the number of hubs to be located.
Constraints (\ref{CCOV5}) are flow conservation constraints, such that the
flow that enters  any hub $k$ with final destination $m$ is the same that the flow
that leaves hub $k$ with destination $m$.
Constraints  (\ref{CCOV4}) ensure that if the final destination site is a hub,
then the flow goes at most through one additional hub.
These constraints are redundant whenever the cost structure satisfies the
triangular inequality, however they are useful in reducing solution times (see \cite{PRRCH11}).
Constraints (\ref{CCOV6}) and (\ref{CCOV7}) establish again that the intermediate nodes in
any origin-destination path should be open hubs.
Constraints \eqref{UUVC9} establish the capacity constraints of the hubs. Observe that this family of constraints make redundant the family (\ref{CCOV2}), but we have kept it because it reduces the computational times.
Constraints \eqref{CCOV8} link sorting and covering variables. They state that the number of allocations with a cost at least $\hat{c}_{(h)}$ must be equal to the number of sites that support shipping costs to the first hub greater than or equal to $\hat{c}_{(h)}$. Finally, constraints \eqref{CCOV9}
are a group of sorting conditions on the variables $\bar{u}_{ih}$.}

The reader may note that this formulation is a natural extension for the capacitated version of the radius formulation already considered  for the uncapacitated ordered median hub location problem in \cite{PRRCH11,PRRCH13}.
However, although this formulation is enough to specify the CSA-OMHLP, we have found that for solving medium sized problems it produces very large MIP models, which are difficult to solve with standard MIP solvers (CPLEX, XPRESS; Gurobi...). Therefore, some alternatives should be investigated.

One way to improve the performance  of the above formulation is to take advantage of some features of that model to reduce the number of variables. In  this case, one can succeed reducing the number of $u$ variables.
The  logic  of the above formulation can be further strengthen for important particular cases of the discrete ordered median hub location problem. In the following, we show a reformulation
that is based on taking advantage of sequences of repetitions in the $\lambda$-vector.
(See \cite{MNV10,PRRCH13, PhDAna}   for  similar reformulations applied to other location problems.)

One can realize that for $\lambda$-vectors with sequences of repetitions --i.e. the center,
$k$-centrum, trimmed means or median among others--, many variables used in formulation
$F_{\bar{u}}$
are not necessary (since they are multiplied by zero in the objective {function}), and some others can be glued together (since they have the same  coefficient in the objective function).
{Moreover, under the assumption of} the  free self-service, and that any origin is allocated to itself if it is a hub, we conclude that the $p$ smallest transportation {costs} from the origin to the first hubs are $0$, i.e. the first $p$ components of the $\lambda$-vector are multiplied  by 0.
Therefore, in order to simplify the problem one can disregard the $p$ first components of the $\lambda$-vector.
Let $\tl =(\tl_1,\ldots,\tl_{n-p}):=(\lambda_{p+1},\ldots,\lambda_n)$.

In order to give a formulation for the CSA-OMHLP taking advantage of these facts, we need to introduce some additional notation.
Let $I$ be the number of
blocks of consecutive  equal non-null elements in $\tl$ and define the vectors:
\begin{enumerate}
\item $\gamma=(\gamma_1,\ldots,\gamma_I)$, being $\gamma_i$, $i=1,\ldots,I$ the value of the
elements in the $i$-th  block of repeated elements in $\tl$.
\item $\alpha=(\alpha_1,\ldots,\alpha_I,\alpha_{I+1})$, being $\alpha_i$ with $i=1,\ldots,I$,
the number of zero entries between the $(i-1)$-th and $i$-th blocks of positive elements in $\tl$
and $\alpha_{I+1}$  the number of zeros, if any, after the $I$-th block of non-null elements in
$\tl$. For notation purposes we define $\alpha_0=0$.
\item $\beta=(\beta_1,\dots,\beta_I)$, being $\beta_i$, $i=1,\ldots,I$ the number of elements in
the $i$-th block of non-null elements in $\tl$. For the sake of compactness, let
$\beta_0=\beta_{I+1}=0$.
\end{enumerate}
Next,  let  denote
$\overline{\alpha}_i=\sum_{j=1}^i \alpha_j$, $\overline{\beta}_i=\sum_{j=1}^i \beta_j$ and
recall that $W_i=\sum_{j \in A} w_{ij}$.
Moreover, for all $i=1,\ldots,I$ and $h=1,\ldots,G$, let us define the following set of decision
variables:
\begin{eqnarray*}
u_{ih}&=&\left\{\begin{array}{ll} 1, & \mbox{if the } \displaystyle (p+\overline{\alpha}_i+
\overline{\beta}_{i-1}+1) \mbox{-th  assignment cost is at least } \hat c_{(h)}, \\
 0, &
\mbox{otherwise.}\end{array} \right.\\
v_{ih}&=&
 \mbox{Number of assignments in the $i$-th block  between the positions } \\
 & &
 \nonumber
 \displaystyle p+ \overline{\alpha}_i+\overline{\beta}_{i-1}+1 \mbox{ and }
 p+\overline{\alpha}_i+\overline{\beta}_{i}
 \mbox{ that are at least } \hat{c}_{(h)}.
\end{eqnarray*}

With the above notation the formulation of CSA-OMHLP is:
\begin{eqnarray}
\label{UVC0}
{(F_{uv})}&\min & \sum_{i=1}^{I} \sum_{h=2}^{G} \gamma_{i}(\hat c_{(h)}-\hat c_{(h-1)})v_{ih}+
\sum_{k \in A} \sum_{\ell \in A} \sum_{m \in A} (\mu c_{k \ell}+\delta c_{\ell m})s_{k\ell m}\\
\nonumber
&s.t. &
Constraints: (\ref{CCOV1})-(\ref{UUVC9}), \\ 
\label{UVC8}
& & \sum_{i=1}^{I} \alpha_ i u_{ih}+ \sum_{i=1}^I v_{ih}+\alpha_{I+1} \geq  \sum_{j \in A}
\sum_{\stackrel{k \in A}{\hat{c}_{jk}  \ge \hat  c_{(h)}}}  x_{jk}, \quad \forall h=2,\ldots,G \\
\label{UVC11}
& & u_{ih}\ge u_{i-1,h}, \quad \forall i=2,\ldots,I, \: h=1,\ldots,G\\
\label{UVC12}
& & \beta_{i-1} u_{ih}\ge v_{i-1,h}, \quad \forall i=2,\ldots,I, \: h=1,\ldots,G \\
\label{UVC13}
& & v_{ih}\ge \beta_{i} u_{ih}, \quad \forall i=1,\ldots, I,\quad h=1,\ldots,G   \\
 \label{UVC23}
& &   u_{ih}\in \{0,1\}, v_{ih}\in \mathbb{Z}\cap [0,\beta_i], \quad \forall  i=1\ldots,I,\quad h=1,\ldots, G \\
\label{UVC24}
& & x_{jk} \in \{0,1\},  s_{k \ell m} \ge 0, \quad \forall j,k, \ell,m\in A.
\end{eqnarray}

The objective function  \eqref{UVC0}  is a reformulation of \eqref{CCOV0} substituting the $\bar{u}$ variables by the new $u$, $v$ variables and the vector $\gamma$, taking advantage of the $\lambda$ vector properties.
Constraints \eqref{UVC8} 
ensure that the number of sites that support a shipping cost to the first hub greater than or equal
to  $\hat{c}_{(h)}$ is either equal to the number of allocations with a cost at least $\hat{c}_{(h)}$
whenever $v_{Ih}>0$ or less than or equal to $\alpha_{I+1}$ otherwise.
Constraints \eqref{UVC11} are sorting constraints on the {variables $u$} similar to constraints \eqref{CCOV9}, and
constraints \eqref{UVC12}-\eqref{UVC13}  provide  upper and lower bounds on  the {variables $v$}
depending on the values of {variables $u$}.

The main difference between $F_{\bar{u}}$ and  $F_{uv}$  is that all $\bar{u}_{ih}$ variables associated with blocks of zero $\lambda$-values are removed, and those associated with each block of non-null $\lambda$ values are replaced by $2\times G$ variables. Therefore, overall we reduce the number of variables by $(n-2I)\times G$.

Note that in Formulation $F_{\bar{u}}$, the family of constraints that links covering variables
(variables $\bar{u}$) and the allocation variables (variables $x$), i.e. \eqref{CCOV8}, is given with equalities. This fact implies that the actual dimension of the feasible region in the space of $\bar{u}$ and $x$ variables is smaller than the one that we were currently  working on. This  is exploited in the new formulation. Indeed, we reduce the number of variables used in the sorting phase replacing $\bar{u}$ by $u$ and  $v$. Therefore, the dimension of the feasible region in the space of $u$, $v$, $x$ variables has smaller dimension. In addition, the constraints that link sorting and design variables, namely \eqref{UVC8}, are given as inequalities. This new representation, although valid for the problem, induces some loss of information in that it does not allow us to take full control of the exact number of {allocations} at some specific cost. This does not affect the resolution process but influences the derivation of valid inequalities.

Finally, for those cases where $\beta_i=1$ we observe that $v_{ih}=u_{ih}$. This set of
constraints whenever valid, was added to reinforce the formulation.

\begin{ejem}\label{PrimerEjemplo}
To illustrate how the $F_{uv}$ versus $F_{\bar{u}}$ formulations work,
we consider the following data. 
Let $A=\{1,\dots,6\}$ be a set of sites and assume that we are
interested in locating $p=2$ hubs. Let the cost and flow matrices
be as follows:

\begin{tiny}
$$
C=\left(\begin{array}{cccccccccccc}
0&14&15&16&15&9\\
5&0&7&2&19&16\\
16&5&0&7&1&19\\
12&1&10&0&13&1\\
1&9&9&15&0&2\\
8&10&16&8&4&0\\
\end{array}\right),\;
W=\left(\begin{array}{cccccccccccc}
0&15&2&8&11&2\\
19&0&1&16&20&7\\
3&9&0&3&11&16\\
7&2&5&0&14&5\\
15&4&20&4&0&1\\
12&4&7&11&18&0\\
\end{array}\right).\;
$$
\end{tiny}
Therefore, $\hat{c}_{(\cdot)}$, the sorted vector of $\hat c$, is in our case \\
$\hat c_{(\cdot)}=
[0, 33, 42, 44, 88, 126, 208, 210, 294, 315, 330, 342, 396, 416, 429, 441, 520, 532, 570, 608, 660,\\  672, 798, 832, 1008, {1197}].$  Hence, $G=26$.
Let
$\lambda=(0,1,0,0,1,1)$, $\mu=0.7$, $\delta=0.9$, and the capacity constraints vector
$b=(119, 119, 113, 145, 149, 140)$.
{The optimal solution opens hubs $4$ and $6$}. The allocation of origin sites to
first hub is given by the following values of the {variables $x$} (see Figure \ref{fig:ex1.1}):
$$x_{16}=x_{24}=x_{34}= x_{44}= x_{56}= x_{66}=1. $$
Analogously, the allocation of first hubs to final destinations are  given by the values of the non null {variables $s$}. Thus,  the flows considering as first hubs $4$ and $6$ are  (see Figure \ref{fig:ex1.1} for a graphical representation of the delivery paths):
\noindent
$$s_{442}=11, \: s_{443}= 6,  \: s_{444}= 19,  \: s_{461}= 29, \: s_{465}=45, \: s_{466}= 28; $$
$$ s_{642}= 23, \:  s_{644}= 23, \: s_{661}= 27, \: s_{663}= 29, \: s_{665}= 29, \: s_{666}= 3 .$$
\begin{figure}[!ht]
\begin{center}
\begin{pspicture}(0,2.5)(8.,5)
\psset{arrows=-,fillstyle=solid,radius=1.5}
\psset{radius=.5,unit=0.75cm}
\cnodeput(-2,2){2b}{\makebox(0, 0.75){$2$}}
\cnodeput(-2,7){1b}{\makebox(0, .75){$1$}}
\cnodeput(4,7){5b}{\makebox(0, .75){$5$}}
\cnodeput(4,2){3b}{\makebox(0, .75){$3$}}
\cnodeput(9,2){2}{\makebox(0, 0.75){$2$}}
\cnodeput(9,7){1}{\makebox(0, .75){$1$}}
\cnodeput(15,7){5}{\makebox(0, .75){$5$}}
\cnodeput(15,2){3}{\makebox(0, .75){$3$}}
\psset{fillstyle = vlines*,hatchcolor=black,hatchsep=8 pt}
\cnodeput(1,3){4b}{\makebox(0, .75){$4$}}
\cnodeput(12,3){4}{\makebox(0, .75){$4$}}
\psset{fillstyle = solid,fillcolor=lightgray!50}
\cnodeput(1,6){6b}{\makebox(0, .75){$6$}}
\cnodeput(12,6){6}{\makebox(0, .75){$6$}}
\psset{fillstyle = none}

\psset{linewidth=1.pt}
\ncline[arrowscale=2]{->}{1b}{6b}
\ncline[arrowscale=2]{->}{2b}{4b}
\nccircle[nodesep=4pt,angle=-180,arrowscale=2]{->}{4b}{.45cm}
\ncline[arrowscale=2]{->}{3b}{4b}
\nccircle[nodesep=4pt,angle=0,arrowscale=2]{->}{6b}{.45cm}
\ncline[arrowscale=2]{->}{5b}{6b}
\psset{linewidth=1.pt,linestyle=dashed,linecolor=black}
\ncarc[arrowscale=2,arcangle=20]{->}{4}{2}
\nccircle[nodesep=4pt,angle=-180,arrowscale=1.]{->}{4}{.35cm}
\psline[]{->}(11.75,3.6)(11.75,5.4)
\psline[arrowscale=2]{->}(11.75,5.4)(9.4,6.5)
\ncarc[arrowscale=2,arcangle=-20]{->}{4}{3}
\psline[arrows=->](12.25,3.6)(12.25,5.4)
\psline[arrows=->,arrowscale=2](12.25,5.4)(14.6,6.5)
\psset{linewidth=1.2pt,linestyle=solid,linecolor=lightgray}
\ncarc[arrowscale=2,arcangle=0]{->}{6}{4}
\ncarc[arrowscale=2,arcangle=20]{->}{6}{5}
\nccircle[nodesep=4pt,angle=0,arrowscale=1]{->}{6}{.35cm}
\ncarc[arrowscale=2,arcangle=-20]{->}{6}{1}
\ncarc[arrowscale=2,arcangle=0]{->}{6}{3}
\psline[]{->}(11.5,5.4)(11.5,3.5)
\psline[arrowscale=2]{->}(11.5,3.5)(9.45,2.5)
\end{pspicture}
\end{center}
\vspace{1cm}
\caption{\footnotesize Illustration of Example \ref{PrimerEjemplo}. Left figure represents the allocations of sites to their corresponding first hubs. Right figure represents the flow pattern to the final destinations from the first hubs:  4 (dashed lines) and  6 (grey lines).
}
\label{fig:ex1.1}
\end{figure}
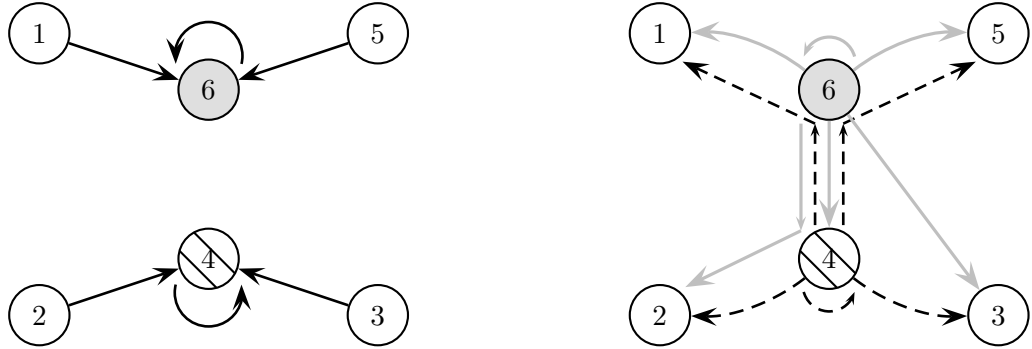

Moreover,  the covering  variables $\bar{u}_{ih}$ are given below. Due to their structure, we only report for each $i$ the  last one and  first zero occurrences since they characterize the remaining values.
$$
\begin{array}{lll}
i=1 \mapsto \bar{u}_{11}= 1, \bar{u}_{12}= 0 &
i=2 \mapsto \bar{u}_{21}= 1, \bar{u}_{22}= 0 &
i=3 \mapsto \bar{u}_{35}= 1, \bar{u}_{36}= 0 \\
i=4 \mapsto \bar{u}_{46}= 1, \bar{u}_{47}= 0 &
i=5 \mapsto \bar{u}_{59}= 1, \bar{u}_{5,10}= 0 &
i=6 \mapsto \bar{u}_{6,12}= 1, \bar{u}_{6,13}= 0. \\
\end{array}$$

The first two assignments are done at a cost $c_{(1)}=0$, corresponding to the two hubs ($p=2$). The next assignment has been done at a cost $c_{(5)}=88$, since $\bar{u}_{35}= 1$ and $\bar{u}_{36}= 0$, and so on. The rest of assignments costs are then $c_{(6)}=126$, $c_{(9)}=294$ and $c_{(12)}=342$.

Hence,  the overall cost of this solution is
$$\displaystyle \sum_{i \in A} \sum_{h=2}^{G} \lambda_{i}(\hat c_{(h)}-\hat c_{(h-1)})\bar{u}_{ih}+
\sum_{k \in A} \sum_{\ell \in A} \sum_{m \in A} (\mu c_{k \ell}+\delta c_{\ell m})s_{k\ell m}=636 +1500.8=2136.8.$$
In addition, to illustrate how the formulation $F_{uv}$ is related with  $F_{\bar{u}}$, we also include the solution of the
covering  variables ${u}_{ih}$ and ${v}_{ih}$:
$$
\begin{array}{l}
I=1 \mapsto {u}_{16}= 1,  u_{17}= 0 \\
I=1 \mapsto {v}_{11}= \ldots = v_{19}=2; \: {v}_{1,10}=\ldots = v_{1,12}=1; \: v_{1,13}=0
\end{array}$$

\begin{figure}[!ht]
$$ \bar{u}_{i,h}=
\begin{array}{c}
  \\
{_1} \\
{_2} \\
{_3} \\
{_4} \\
{_5} \\
{_6} \\
\end{array}
\begin{array}{c}
 {\tiny \begin{array}{cccccccccc}
  c_{(1)} &  c_{(2)} &  c_{(3)}& c_{(4)}&c_{(5)}& c_{(6)}&\ldots &  c_{(24)}&  c_{(25)}&  c_{(26)}
 \end{array}} \\
\left(
\begin{array}{cccccccccc}
\hspace*{0.15cm} 1\hspace*{0.15cm}  & \hspace*{0.15cm}0& \hspace*{0.15cm}0& \hspace*{0.15cm}0& \hspace*{0.15cm}0 \hspace*{0.15cm}& \hspace*{0.15cm}0\hspace*{0.15cm}&\hspace*{0.15cm} \ldots \hspace*{0.15cm}&\hspace*{0.15cm}0\hspace*{0.15cm}&\hspace*{0.15cm}0\hspace*{0.15cm}&\hspace*{0.15cm} 0 \hspace*{0.15cm}\\
1 & 0 & 0& 0& 0& 0& \ldots &0&0& 0 \\
1 & 1& 1& 1& 1 & 0& \ldots &0&0& 0 \\
1 & 1& 1& 1& 1 & 1& \ldots &0&0& 0 \\
1 & 1& 1& 1& 1 & 1& \ldots &0&0& 0 \\
1 & 1& 1& 1& 1 & 1& \ldots &0&0& 0 \\
\end{array}\right)\\
\end{array}
\Rightarrow
\lambda=
\begin{array}{c}
\begin{array}{lrl}
& &
\end{array}\\
\left(
\begin{array}{c}
0 \\
1 \\
0 \\
0 \\
1 \\
1
\end{array}
\right)
\end{array}
\begin{array}{c}
\begin{array}{lrl}
 & &
\end{array}\\
\hspace*{-0.5cm}
\begin{array}{ll}
\left\} \hspace*{-2cm} \begin{array}{rl} {\it }&\\ {\it }\hspace*{-2cm}\end{array}\right. &   p  \\
\left\}\hspace*{-2cm} \begin{array}{rl} { } &\\ { }\hspace*{1cm}\end{array}\right. &\alpha_1   \\
\left\}\hspace*{-2cm} \begin{array}{rl} { }& \\ { }\end{array}\right. & \beta_1  \left\}\begin{array}{r} v_1\end{array}
\right.
\end{array}
\end{array}
\begin{array}{c}
\begin{array}{lrl}
 & &
\end{array}\\
\hspace*{-2.9cm}
\begin{array}{ll}
&  \\
&  \\
&  \\
&  \\
\longrightarrow\longrightarrow u_1 & \\
&  \\
\end{array}
\end{array}
$$
\caption{Variables and lambda vector of Example \ref{PrimerEjemplo}\label{t:ejemplo}}
\end{figure}

Note that we have only one block of repeated non-null elements of the $\tl$-vector, so $I=1$. (See the right part of Figure \ref{t:ejemplo}.)
The number of zero entries between  two blocks is $\alpha_1=2$, and the number of elements in the $1$-st block  of non-null elements is $\beta_1=2$.
 Furthermore, $\gamma_1=1$ is the repeated value in the $1$-st block.
 \vspace*{-0.5cm}
$$
\begin{array}{c}
{\tiny
\begin{array}{p{.5cm}p{.5cm}p{.5cm}p{.5cm}p{.5cm}p{.5cm}p{.5cm}p{.5cm}p{.5cm}p{.5cm}p{.5cm}p{.5cm}p{.5cm}p{.5cm}p{.5cm}p{.5cm}p{.5cm}p{.5cm}}
&$c_{(1)}$ & $c_{(2)}$ &$c_{(3)}$ & $\ldots$ &$c_{(9)}$&$c_{(10)}$&$c_{(11)}$&$c_{(12)}$&$c_{(13)}$&$c_{(14)}$&$c_{(15)}$&$c_{(16)}$&$c_{(17)}$& $\ldots$ &$c_{(25)}$&$c_{(26)}$
\end{array}
}
\\
{u}_{1,h}=(
\begin{array}{p{.5cm}p{.5cm}p{.5cm}p{.5cm}p{.5cm}p{.5cm}p{.5cm}p{.5cm}p{.5cm}p{.5cm}p{.5cm}p{.5cm}p{.5cm}p{.5cm}p{.5cm}p{.5cm}p{.3cm}}
$1$& $1$&$1$& \ldots &$1$&$1$&$1$&$1$&$1$&$1$&$1$&$1$&$0$& \ldots &$0$&$0$
\end{array} )\\
{v}_{1,h}=(
\begin{array}{p{.5cm}p{.5cm}p{.5cm}p{.5cm}p{.5cm}p{.5cm}p{.5cm}p{.5cm}p{.5cm}p{.5cm}p{.5cm}p{.5cm}p{.5cm}p{.5cm}p{.5cm}p{.5cm}p{.3cm}}
$2$& $2$&$2$& \ldots &$2$&$1$&$1$&$1$&$0$&$0$&$0$&$0$&$0$& \ldots &$0$&$0$
\end{array} )\\
\end{array}
$$
The variable ${u}_{1,h}$ points out to the row  $p+\alpha_1+1 = 5$ of the original variable $\bar{u}_{i,h}$.
Whereas the variable ${v}_{1,h}$ accounts for the number of assignments between the positions
 $p+\alpha_1 +1 = 5$  and $p+\alpha_1+\beta_1 = 6$ of  $\bar{u}_{i,h}$ that are at least $\hat{c}_{(h)}$. (See Figure  \ref{t:ejemplo}.)

From $\bar{u}_{i,h}$, we know that the $5$-th assignment cost is $\hat{c}_{(9)}$ and the $6$-th assignment cost is $\hat{c}_{(12)}$. For this reason ${v}_{1,h}$ is equal to $2$ up to column $9$, this is the number of assignment costs greater  than $\hat{c}_{(9)}$. Being this number equal to 1 from $h=10$ to $h=12$, and zero for the remaining columns.

 Applying this formulation, the overall reduction in the number of variables is $(n-2I)\times G=104$.
The rest of variables $x$ and $s$ remain the same, and again the overall cost of this solution is
$$\displaystyle \sum_{i=1}^{I} \sum_{h=2}^{G} \gamma_{i}(\hat c_{(h)}-\hat c_{(h-1)})v_{ih}+
\sum_{k \in A} \sum_{\ell \in A} \sum_{m \in A} (\mu c_{k \ell}+\delta c_{\ell m})s_{k\ell m}={636 +1500.8=2136.8}.$$

\end{ejem}

\section{Strengthening the formulation}\label{Sec:StrenForm}

\subsection{Variable fixing}\label{SubSec:StrenForm1}
\label{ssec:vf}
Next, we describe some preprocessing procedures that we have applied to reduce further the
size of  formulation $F_{uv}$. We present a number of variable fixing possibilities for the set of
 {variables $u$ and $v$}  which are useful in the overall solution process.
The variable fixing procedures developed in this section  are based on ideas used in \cite{PRRCH11,PRRCH13} and taking advantage of the capacity constraints.  Indeed, we are adding the reinforced effective capacity constraints, { as well as some surrogated version of constraints} \eqref{PREC2} since
in this aggregated form they give better running times.
The  preprocessing phase developed in this paper also provides new upper
and lower bounds on the $v$ variables.
The percentage of variable reduction obtained by these procedures can be found in Tables \ref{tabla1} and \ref{tabla2} (column named as 'Fixed').

Before describing these procedures for fixing variables,   the following simple arguments allows us to fix some variables:
\begin{enumerate}
\item First, $c_{jj}=0$ $\forall j\in A$, i.e., $\hat{c}_{(1)}=0$. Moreover,  every origin where it has been located a hub will be allocated to itself as a first hub.
\item Second,  $\hat{c}_{jk}\ne 0$  if and only if   $ j \ne k$, i.e., any non-hub origin is allocated to a first hub at  a cost of at least   $\hat{c}_{(2)}$.
\end{enumerate}
Therefore, since in this formulation the first $p$-allocations are considered only implicitly, we can fix  $u_{i1}=1$,  $v_{i1}=\beta_i$ as well as $u_{i2}=1$,
$v_{i2}=\beta_i$, $ \forall i=1,\ldots,I.$

\subsubsection{Preprocessing Phase 1:  { Fixing variables to the upper} bounds for the formulation with covering variables strengthen with capacity constraints.}
\label{ssec:vf1}

Due to the definition of the variables in formulation $F_{uv}$, one can expect that  $u_{ih}=1$ whenever
$i$ is large and $h$ is small to medium size because this would mean that the
$(p+\overline{\alpha}_{i}+\overline{\beta}_{i-1}+1)$-th sorted allocation cost would not have been done at cost less than $\hat{c}_{(h)}$.
The reader should observe that an analogous strategy applies to the {variables $v$} since their interpretation is similar, but
in this case the values of $v_{ih}$ would be fixed to $\beta_i$.
For the cases where it is not possible to fix the corresponding $v$-variable, it could be still possible
to establish some  lower bounds as we will see later.

Next, to fix {variables $u_{ih}$ and $v_{ih}$} for $i=1,\ldots,I$, $h =1,\ldots,G$,  we deal with
an auxiliary problem that maximizes the number of variables that may assume zero values, satisfying
$\hat{c}_{jk} \le \hat c_{(h-1)}$ and the capacity constraints.
For any $h =1,\ldots,G$ and $j,k\in A$  such that  $\hat c_{jk}\le \hat c_{(h-1)} $ let
\begin{equation}
\label{zdef}
z_{jk}^h=
\begin{cases}
1, & \mbox{if origin site $j$ is assigned to hub $k$} \\
0, & \mbox{otherwise.}
\end{cases}
\end{equation}
To avoid possible misunderstanding in the cases where variables  $z_{jk}^h$ are not defined,
i.e. when $\hat c_{jk}>\hat c_{(h-1)}$,  we can assume that $z_{jk}^h:=0$.

For a given $h$, we introduce the {\it effective capacity} of a hub $k$ at a cost at most $\hat c_{(h-1)}$ as,
\begin{equation}\label{eqEffCap}
\displaystyle {b_k^{h-1}}:= \min\left\{b_k, \sum_{ \stackrel{s \in A}{\hat c_{sk}\le \hat c_{(h-1)}}} W_s
\right\} .
\end{equation}
Indeed, the capacity of a hub $k$ is always lower than or equal to $b_k$. In addition, if we restrict ourselves to the nodes served at a cost of at most $\hat c_{(h-1)}$, then the actual capacity to cover this set should be lower than $ \sum_{ s \in A:\hat c_{sk}\le \hat c_{(h-1)}} W_s$,  and this gives us the expression of ${b_k^{h-1}}$.

The optimal value $P_1(h)$ of the following problem fixes the maximal number of  allocations that may be feasible at a cost of at most $\hat c_{(h-1)}$.
\begin{eqnarray}
\label{PREC0}
P_1(h):=&\max &{\sum_{\stackrel{j,k \in A:}{\hat c_{jk}\le \hat c_{(h-1)}}}} z^h_{jk} \nonumber\\
&s.t. &
\label{PREC1} {\sum_{{k \in A:}{\hat c_{jk}\le \hat c_{(h-1)}}}}
z^h_{jk}\leq1, \quad \forall j\in A,  \nonumber\\
\label{PREC2}
& &{\sum_{{j \in A:}{\hat c_{jk}\le \hat c_{(h-1)}}}}
z^h_{jk} \leq ny_{k}, \quad \forall k\in A \\
\label{PREC3}
& &\sum_{k \in A} y_{k}\le p \nonumber\\
\label{PREC4}
& &  {\sum_{j \in A:\hat c_{jk}\le \hat c_{(h-1)}}}
W_jz^h_{jk} \le {b_k^{h-1}}y_{k}, \quad \forall k\in A \nonumber\\
\label{PREC5}
& &  z^h_{jk}, y_k \in \{0,1\},  \quad \forall j,k\in A ,\; \forall h =1, \ldots ,G. 
\nonumber
\end{eqnarray}

Then, depending on the value
$P_1(h)$ we can  fix  some variables to their upper bounds. Let us denote by $i_1(h) \in \{1,\ldots,I\}$ the index such that
$$ p+\overline{\alpha}_{i_1(h)-1}+\overline{\beta}_{i_1(h)-1}
<P_1(h) \le p+\overline{\alpha}_{i_1(h)}+\overline{\beta}_{i_1(h)}.$$
$$\left\{\begin{array}{ll}
 u_{ih}=1, v_{ih}=\beta_i, \; i=i_1(h),\ldots,I, &
\mbox{if } \displaystyle P_1(h) \le p+\overline{\alpha}_{i_1(h)}+\overline{\beta}_{i_1(h)-1}  \\
\left\{ \begin{array}{l}
\displaystyle  v_{i_1(h), h}\ge p+\overline{\alpha}_{i_1(h)}+\overline{\beta}_{i_1(h)} -P_1(h),  \\
   u_{ih}=1,\quad v_{ih}=\beta_i, \; i=i_1(h)+1,\ldots,I,
   \end{array} \right\}& \mbox{otherwise.}
\end{array}
\right.
$$

\subsubsection{Preprocessing Phase 2:   Fixing variables to their lower bounds}
\label{ssec:vf2}

Following similar argument to the previous subsection, one can expect that many
 {variables $u$ and $v$}  in the top-right hand corner of the matrices of  {variables $u$ and $v$} , respectively, will take
value 0 in the  optimal solution. Indeed, $u_{ih}=0$ means that the
$(p+\overline{\alpha}_{i}+\overline{\beta}_{i-1}+1)$-th sorted allocation cost is less than
$\hat{c}_{(h)}$ which is very likely to be true if $h$ is sufficiently large and $i$ is small.  Note that an analogous strategy,  applies to the {variables $v$} since their interpretation is
similar. For the cases where it is not possible to fix the corresponding $v$-variable it could be still possible  to establish some  upper bounds.

For any $h=2,\ldots,G$, $j,k\in A$  such that  $\hat c_{jk}\ge \hat c_{(h-1)} $  let define
variables $z_{jk}^h$ as \eqref{zdef}.
To avoid possible misunderstanding in the cases where  variables $z_{jk}^k$ are not defined, i.e. when  ${\hat c_{jk}<\hat c_{(h-1)}}$, we can assume that ${z_{jk}^h}:=0$.
Using these variables, the formulation of the problem that maximizes the number of non-fixed allocations at a cost at most
$\hat c_{(h-1)}$ is:
\begin{eqnarray}
P_2(h):=&\max & {\sum_{\stackrel{j,k \in A:}{\hat c_{jk}\ge \hat c_{(h-1)}}}}  {z_{jk}^h} \nonumber\\
&s.t.&
{\sum_{{k \in A:}{\hat c_{jk}\ge \hat c_{(h-1)}}}} {z_{jk}^h} \le 1, \quad \forall j\in A \nonumber\\
& & {\sum_{{j \in A:}{\hat c_{jk}\ge \hat c_{(h-1)}}}} {z_{jk}^h} \le n y_k, \quad \forall j,k\in A \label{PREP2}\\
& & \sum_{k \in A}  y_k \le p, \nonumber\\
& &   {z_{jk}^h}, y_k \in \{0,1\},\quad  \forall j,k\in A,\; \forall h =1, \ldots ,G.
\nonumber
\end{eqnarray}

Note that the value $P_2(h)$ implies that there are no feasible solutions of the original problem with less than $n-P_2(h)$ allocations fixed at a cost at most $\hat c_{(h)}$.

Let $1 \le i_2(h)\le I$ be the index such that
$$ p+\overline{\alpha}_{i_2(h)-1} +\overline{\beta}_{i_2(h)-1} < n-P_2(h) \le p+\overline{\alpha}_{i_2(h)}+\overline{\beta}_{i_2(h)}.$$
Thus, in any feasible solution of the problem we have that:
$$\left\{\begin{array}{ll}
u_{ih}=0, v_{ih}=0, \: i=1,\ldots,i_2(h) -1, &
\mbox{if }
\displaystyle n-P_2(h) \le p+\overline{\alpha}_{i_2(h)}+\overline{\beta}_{i_2(h)-1} \\
\left\{ \begin{array}{l}
\displaystyle
u_{i_2(h),h}=0, \quad v_{i_2(h),h}\le p+\overline{\alpha}_{i_2(h)}+\overline{\beta}_{i_2(h)}-(n-P_2(h)),  \\
u_{ih}=0, \quad v_{ih}=0,\quad i=1,\dots,i_2(h) -1,
\end{array} \right\} & \mbox{otherwise.}
\end{array}
\right.
$$
Note that whenever $n-P_2(h)=p$ then there is nothing to fix and therefore no variables are set to zero in  column $h$.

\subsection{Valid Inequalities}\label{SubSec:StrenForm2}

In order to strengthen formulation $F_{uv}$ we have studied several families of
valid  inequalities.
In fact, taking advantage of previous experience on the non-capacitated version of the problem we have borrowed 
a first family of valid inequalities that are very simple and that have proven to be effective in different ordered  median problems with covering variables
 \cite{PRRCH11,PRRCH13}.
This family is
\begin{eqnarray}\label{in:nuv}
u_{ih} \ge u_{i, h+1}, & i=1,\ldots, I,\; h=1,\ldots,G-1, \\
v_{ih} \ge v_{i, h+1}, & i=1,\ldots, I,\; h=1,\ldots,G-1. \label{in:nuv2}
\end{eqnarray}
Since, these families are straightforward consequence of the definition of  {variables $u$ and $v$}  they have been included in the original formulation.

In the following,  we describe several alternative families of valid inequalities: three sets of inequalities, \eqref{in:DV3_3},\eqref{in:DV3_1}-\eqref{in:DV3bis_2}, and \eqref{in:DV3b_1}-\eqref{in:DV3c_2}, based on the combination of ordering and capacity requirements and two more sets, \eqref{in:cortepacking1} and \eqref{const1}-\eqref{const3}, that do not use capacities.

\subsubsection{First family of valid inequalities: Valid inequalities based on capacity I }

We can add to this model several  families of valid inequalities based on capacity issues that help in solving the problem by reducing the gap of the linear relaxation and the CPU time to explore the branch and bound search tree.

Observe that the  capacity of the set of hubs that may be used to assign origins in $A$ at a cost at most $\hat c_{(h-1)}$, is given by
$$ \displaystyle \sum_{k \in A}{b_k^{h-1}}x_{kk}$$
where $\displaystyle {b_k^{h-1}}$ is the effective capacity at a cost at most $\hat c_{(h-1)}$, defined by \eqref{eqEffCap}.
Recall that, although the capacity of a hub $k$ is always lower than or equal to $b_k$, when we restrict to the nodes served at a cost of at most $\hat c_{(h-1)}$, then the actual capacity to cover this set should be lower than $b_k^{h-1}$.
 Making use of the above observation, we can add the following family of constrains as valid inequalities
\begin{eqnarray}
\label{in:DV3_3}
\sum_{\stackrel{j \in A}{\hat{c}_{jk} \le \hat{c}_{(h-1)}}}  W_j x_{jk} \le  {b_k^{h-1}} x_{kk} \quad \forall h=2,\ldots, G ,  k\in A
 \end{eqnarray}
which enforces that all the flow sent from origin-hubs at a cost at most $\hat c_{(h-1)}$ cannot exceed
the effective capacity at that cost. Observe that in the case where $b_k^{h-1}$ takes the value $b_k$ (11) dominates {\eqref{in:DV3_3}}, but in the case where
$b_k^{h-1} = \sum_{ {s \in A:}{\hat c_{sk}\le \hat c_{(h-1)}}} W_s$, {\eqref{in:DV3_3}} becomes
$\sum_{{j \in A:}{\hat{c}_{jk} \le \hat{c}_{(h-1)}}} W_j x_{jk} \le  \sum_{{j \in A:}{\hat{c}_{jk} \le \hat{c}_{(h-1)}}} W_j x_{kk}$.
Observe that this last valid inequality is an alternative surrogation,  with capacity coefficients, of constraints
$x_{jk}\le x_{kk}$ that although valid do not appear in the model because of its large cardinality. This new form of
aggregation has provided good results in the computational experiments.

\subsubsection{Second family of valid inequalities: Valid inequalities based on capacity II }

This section introduces another family of valid inequalities based on capacity issues that help in solving  the problem.
In order to  present these new valid inequalities based on capacity requirements we introduce some
new notation. Assume without
loss of  generality that  $W_i \le W_{i+1}$  for $i\in A\setminus \{n\}$ and let
$\overline{W}_j=\sum_{r=1}^j W_r$ and $S_k:=\{i\in A:   i\le k\}$ for $k\in A$ be a given set of
origin sites.

In case that the effective capacity at a cost at most $\hat c_{(h-1)}$  is not sufficient to cover
the  demand of  $S_k$, i.e. $ \displaystyle  \sum^{n}_{j=1}  \Big(b_j^{h-1}-W_j\Big) x_{jj}$ is less
than  $\sum_{s=1}^k  W_s(1-x_{ss})$, then  at most $k-1$ origins  of $S_k$ can  be allocated at a cost lower
than or equal to $\hat c_{(h-1)}$. This  argument  can be applied for each $h$ to the corresponding
$\hat c_{(h-1)}$ value.
Moreover, we have chosen this particular structure of $S_k$ consisting of the $k$ origins (nodes) with  the $k$-smallest flows, because   given a fixed amount of flow, the maximal cardinality set of origins, such that the overall flow originated in this set is lower than or equal to this amount, is provided by a  set $S_k$ for some $k\in A$. Therefore, since we are dealing with the worst cases, it allows us to fix some {variables $u$ and  $v$} through the following valid inequalities.
For each $h=2,\ldots, G$ we obtain the following:
\begin{itemize}
\item
If   $\overline{\alpha}_{i-1}+\overline{\beta}_{i-1} < k \le \overline{\alpha}_i+\overline{\beta}_{i-1}+1$,   for some $i \in \{1,\ldots , I\}$, 
 namely if the index of the last element, $k$, that defines $S_k$ lies in the $i$-th block of null elements  in the $\tilde \lambda$ vector, then
\begin{equation} \label{in:DV3_1}
\overline{W}_k u_{ih}+\sum_{j \in A}\Big(b_j^{h-1}-W_j\Big)x_{jj} \ge \sum_{s=1}^k  W_s(1-x_{ss}).
\end{equation}
Observe that the above inequality amounts to a disjunctive condition: either the effective capacity at a cost at most $\hat c_{(h-1)}$ is enough to cover the demand of the $k$ smallest flows from origin sites or the $i$-th sorted cost allocation must be assigned at a cost at least $c_{(h)}$.
\item
If $\overline{\alpha}_i+\overline{\beta}_{i-1}+1 < k \le \overline{\alpha}_i+\overline{\beta}_{i}$, and 
namely if the index of the last element, $k$, that defines $S_k$ lies in the $i$-th block of non-null elements, then
\begin{equation} \label{in:DV3_2}
\overline{W}_k v_{ih}+\sum_{j \in A}\Big(b_j^{h-1}-W_j\Big)x_{jj} \ge \sum_{s=1}^k W_s(1-x_{ss}).
\end{equation}
In this case, the inequality is similar to the previous one but written in terms of the variables $v$ that allow to control the capacity whenever $k$ falls within a block of non-null elements in the $\tilde \lambda$ vector.
\end{itemize}
\begin{rem}
Recall that  if $k > P_1(h)-p$, variables $u_{ih}$ and $v_{ih}$ have been already fixed by the Preprocessing Phase 1, and for the above inequality to be effective
$k \le P_1(h)-p$, or equivalently, $i\le i_1(h)$ .
\end{rem}

Based on the same  arguments  we can add a larger family of valid inequalities built on arbitrary sets of origin sites. Let $S$
be a set of origin sites, and suppose that $A_S=\sum_{s \in S} W_s$ satisfies $\overline W_k\le A_S < \overline  W_{k+1}$ ,
\begin{itemize}
\item If
$\overline{\alpha}_{i-1}+\overline{\beta}_{i-1} < k \le \overline{\alpha}_i+\overline{\beta}_{i-1}+1$, for some $i \in \{1,\ldots , I\}$, and $k \le P_1(h)-p$
then
\begin{equation} \label{in:DV3bis_1}
A_S u_{ih}+\sum_{j \in A}\Big(b_j^{h-1}-W_j\Big)x_{jj} \ge \sum_{s=1}^k  W_s(1-x_{ss}).
\end{equation}
\item
If $\overline{\alpha}_i+\overline{\beta}_{i-1}+1 < k \le \overline{\alpha}_i+\overline{\beta}_{i}$, and $k \le P_1(h)-p$ then
\begin{equation} \label{in:DV3bis_2}
A_S v_{ih}+\sum_{j \in A}\Big(b_j^{h-1}-W_j\Big)x_{jj} \ge \sum_{s=1}^k W_s(1-x_{ss}).
\end{equation}
\end{itemize}

Now, assuming a more general case and for an improvement of the above inequalities (\ref{in:DV3_1}) and (\ref{in:DV3_2}),
for a given $h \in \{2,\ldots, G\}$, and  $k \le P_1(h)-p$,
 $k \in \{1,\ldots,\overline{\alpha}_{I+1}+\overline{\beta}_{I})\}$,
let
$$M_s:=\min_{j (\ne s) \in A} \hat c_{js}.$$
$M_s$ is the minimum allocation cost to $s$ as an open hub. In other words, no allocation to hub $s$ is possible at a cost less than $M_s$,
except in case $s$ were a hub itself. We shall call this value the { \it empty radius of $s$}.

Define  $s(h-1,k)$ to be the index of the sorted sequence of elements $W_s$ such that there are exactly $k$ elements $W_s$ with $s\le s(h-1,k)$  and $ M_s \le \hat c_{(h-1)}$, namely $s(h-1,k)$, is the index such that
$$|\{ s \: : \:   s \le s(h-1,k); \: M_s \le \hat c_{(h-1)} \}|=k. $$
Then it holds that,
\begin{itemize}
\item
If    $\overline{\alpha}_{i-1}+\overline{\beta}_{i-1} < k \le \overline{\alpha}_i+\overline{\beta}_{i-1}+1$, for some $i \in \{1,\ldots , I\}$, and  $k \le P_1(h)-p$
\begin{eqnarray} \label{in:DV3b_1}
\left(\overline{W}_{s(h-1,k)}-\!\!\!\!\!\!\!\!\sum_{\stackrel{s=1}{M_s >\hat c_{(h-1)}}}^{ s(h-1,k) } \!\!\!\!\!\!\!\! W_s \right) u_{ih}+	 \sum_{j \in A}\Big(b_j^{h-1}-W_j\Big)x_{jj} \ge  \nonumber \\
\!\!\!\!\!\!\!\!\sum_{\stackrel{s=1}{M_s \le \hat c_{(h-1)}}}^{s(h-1,k) }  \!\!\!\!\!\!\!W_s (1-x_{ss})-
\!\!\!\!\!\!\!\!\sum_{\stackrel{s=1}{M_s > \hat c_{(h-1)} }}^{s(h-1,k)} \!\!\!\!\!\!\!\! { (W_{s(h-1,k)} - W_s)}x_{ss}.
\end{eqnarray}
The above inequality is also a disjunctive condition that ``reinforces''  the  family of valid inequalities (\ref{in:DV3_1}).
It states that if the effective capacity at a cost at most $\hat c_{(h-1)}$ is not enough to
cover the flow sent from origin sites that are not hubs and  that can be allocated at some costs less than or equal to $\hat c_{(h-1)}$ then some of the origin sites with allocation costs less than $\hat c_{(h-1)}$ must be assigned at a cost at least $\hat c_{(h)}$. We observe that the use of $u$ variables in the left-hand side of the inequality is due to the fact that  $k$ falls within a block of null elements in the $\tilde \lambda$ vector. A similar inequality also holds when $k$ falls within a block of non-null elements as shown below.

\item If $\overline{\alpha}_i+\overline{\beta}_{i-1}+1 <k \le \overline{\alpha}_i+\overline{\beta}_{i}$, and $k \le P_1(h)-p$ then
\begin{eqnarray} \label{in:DV3b_2}
\left(\overline{W}_{s(h-1,k)}-\!\!\!\!\!\!\!\!\sum_{\stackrel{s=1}{M_s >\hat c_{(h-1)}}}^{ s(h-1,k) } \!\!\!\!\!\!\!\! W_s \right) v_{ih}+	 \sum_{j \in A}\Big(b_j^{h-1}-W_j\Big)x_{jj} \ge  \nonumber \\
\!\!\!\!\!\!\!\!\sum_{\stackrel{s=1}{M_s \le \hat c_{(h-1)}}}^{s(h-1,k) }  \!\!\!\!\!\!\!W_s (1-x_{ss})-
\!\!\!\!\!\!\!\!\sum_{\stackrel{s=1}{M_s > \hat c_{(h-1)} }}^{s(h-1,k)} \!\!\!\!\!\!\!\! (W_{s(h-1,k)} - W_s) x_{ss}.
\end{eqnarray}
This inequality is similar to the previous one whenever the index $k$ falls within a block of non-null elements in the $\tilde \lambda$ vector.
\end{itemize}

Finally, as the index $s(h-1,k)$ should be greater than or equal to $k$, we can split the above equations, (\ref{in:DV3b_1}) and (\ref{in:DV3b_2}), into $s(h-1,k) - k$ equivalent inequalities.
This is, for any $t = k,\ldots , s(h-1,k)$, define $\hat s(h-1,k,t)$ to be the index
of the sorted sequence of elements $W_s$ such that
 $$|\{ s \:  \:  ; s \le \hat s(h-1,k,t); \: M_s \le \hat c_{(h-1)} \}|+t-\hat s(h-1,k,t)=k.$$
Then it holds that,
 \begin{itemize}
\item
If    $\overline{\alpha}_{i-1}+\overline{\beta}_{i-1} < k \le 	 \overline{\alpha}_i+\overline{\beta}_{i-1}+1$, for some $i \in \{1,\ldots , I\}$, and $k \le P_1(h)-p$
\begin{eqnarray} \label{in:DV3c_1}
\left(\overline{W}_t \,\,\ -\!\!\!\!\!\!\!\!\sum_{\stackrel{s=1}{M_s >\hat c_{(h-1)} \mbox{\scriptsize and } s \le \hat s(h-1,k,t)}}^t\!\!\!\!\!\!\!\! W_s \right) u_{ih}+	 \sum_{j \in A}\Big(b_j^{h-1}-W_j\Big)x_{jj} \ge \nonumber \\
\!\!\!\!\!\!\!\!\sum_{\stackrel{s=1}{M_s \le \hat c_{(h-1)} \mbox{\scriptsize or } s > \hat s(h-1,k,t)} }^t \!\!\!\!\!\!\!\!\!\!\!\!\!\!\!\!W_s (1-x_{ss})
\,\,\, -
\!\!\!\!\!\!\!\!\sum_{\stackrel{s=1}{M_s > \hat c_{(h-1)} \mbox{\scriptsize and } s \le \hat s(h-1,k,t)} }^t \!\!\!\!\!\!\!\!\!\!\!\!\!\!\!\! (W_t - W_s)x_{ss}.
\end{eqnarray}
\item If $\overline{\alpha}_i+\overline{\beta}_{i-1}+1 < k \le \overline{\alpha}_i+\overline{\beta}_{i}$ and $k \le P_1(h)-p$, then
\begin{eqnarray} \label{in:DV3c_2}
\left(\overline{W}_t \,\,\ -\!\!\!\!\!\!\!\!\sum_{\stackrel{s=1}{M_s >\hat c_{(h-1)} , s \le \hat s(h-1,k,t)}}^t\!\!\!\!\!\!\!\! W_s \right)v_{ih}+	 \sum_{j \in A}\Big(b_j^{h-1}-W_j\Big)x_{jj} \ge  \nonumber \\
\!\!\!\!\!\!\!\!\sum_{\stackrel{s=1}{M_s \le \hat c_{(h-1)} \mbox{\scriptsize or } s > \hat s(h-1,k,t)} }^t \!\!\!\!\!\!\!\!\!\!\!\!\!\!\!\!W_s (1-x_{ss})
\,\,\, -
\!\!\!\!\!\!\!\!\sum_{\stackrel{s=1}{M_s > \hat c_{(h-1)} , s \le \hat s(h-1,k,t)} }^t \!\!\!\!\!\!\!\!\!\!\!\!\!\!\!\!(W_t - W_s)x_{ss}.
\end{eqnarray}
\end{itemize}

Observe that if $t=s(h-1,k)$ then $\hat s(h-1,k,t)=s(h-1,k)$. Thus,
the families of valid inequalities \eqref{in:DV3c_1} and \eqref{in:DV3c_2} include as particular instances the families \eqref{in:DV3b_1} and \eqref{in:DV3b_2}.

\subsubsection{Third family of valid inequalities: Disjunctive implications}
The third  family of valid inequalities,  directly borrowed from \cite{PRRCH13}, state disjunctive implications on the origin-first hub allocation costs.
They ensure that either origin site $j$ is allocated to a first hub at a cost of at least $\hat{c}_{(h)}$ or there
is an open hub $k$ such that $\hat{c}_{jk} < \hat{c}_{(h)}$. This  argument  can be formulated through the following family of valid inequalities:
\begin{equation} \label{in:cortepacking1}
\sum_{k \in A: \hat{c}_{jk} \ge \hat{c}_{(h)}} x_{jk}+\sum_{k\in A: \, \hat{c}_{jk} < \hat{c}_{(h)}} x_{kk}
\ge 1, \quad \forall j\in A, \; h=1,\ldots, G.
\end{equation}

\subsubsection{Fourth family of valid inequalities}

Using the definition of the {variables $u$ and $v$}, we establish a lower and an upper  bound of the number of feasible allocations  at a cost $\hat{c}_{(h-1)}$. Observe that using the family of constraints
\eqref{CCOV8} for the original formulation, the exact number of allocations done at a cost
$\hat{c}_{(h-1)}$ is given by { $\sum_{i \in A} (\bar{u}_{i,h-1}-\bar{u}_{i,h})$}. However, since in formulation $F_{uv}$ the number of variables has been considerably reduced, some information is lost. In particular, we cannot keep under control with this new formulation the exact number of allocations at a cost $\hat{c}_{(h-1)}$. Indeed, the counterpart to equalities
\eqref{CCOV8} in formulation $F_{uv}$ is the family of constraints \eqref{UVC8}.
Therefore, we are only able to give a lower and upper bound of this number of allocations. These lower and upper bounds are formulated, respectively, by the following two families of constraints:
\begin{eqnarray}
\label{const1}
\sum_{j \in A} \sum_{\stackrel{k \in A}{\hat c_{jk}= \hat c_{(h-1)}}}
x_{jk} & \geq &
\sum_{i=1}^{I} (v_{i,h-1} - v_{ih}) + \sum_{i=2}^{I} \alpha_i(u_{i-1,h-1} - u_{ih})  , \quad \forall \, h= 2,\ldots , G, \\
\sum_{j \in A} \sum_{\stackrel{k \in A}{\hat c_{jk}= \hat c_{(h-1)}}} x_{jk} & \leq &
\nonumber
\sum_{i=1}^{I} (v_{i,h-1} - v_{i,h}) + \sum_{i=1}^{I} \alpha_i(u_{i,h-1} - u_{i,h})+
(1-u_{Ih})\alpha_{I+1} +  \\
\label{const3}
 & & \alpha_{1}u_{1h} +  \sum_{i=1}^{I-1} \alpha_{i+1}(u_{i+1,h} - u_{ih})
, \qquad \qquad \forall \, h= 2,\ldots , G.
\end{eqnarray}

The first sum in the right hand side of both families gives the exact number of allocations at cost
$\hat{c}_{(h-1)}$ in the positions corresponding to non-null blocks of the vector $\lambda$.
However, the second sum in the right hand side of constraints \eqref{const1} provides a lower
bound of the number of allocations at cost $\hat{c}_{(h-1)}$ in the positions corresponding to the null blocks. In the same way, the $2$nd to the $4$-th sums in  \eqref{const3} provide an upper bound on the number of these allocations.

\section{Computational Results}\label{Sec:Results}

\label{sec:bb_results}

The formulations given to the CSA-OMHLP with the corresponding  strengthening and preproccesing phases, described in this paper, were implemented in the commercial solver XPRESS
IVE 1.23.02.64 running on a
Intel(R) Core(TM) i5-3450 CPU  @3.10GHz  6GB RAM.

The cut generation option of XPRESS was disabled in order to compare the relative performance of the formulations cleanly.

For this purpose we use the AP
data set publicly available at http://www.cmis.csiro.au/or /hubLocation (see \cite{EK96}).
As in previous papers on the field related to the uncapacitated version of this problem, we tested the formulations on a testbed of five
instances for each combination of   costs matrices varying: (i) $n$ in $\{15,20,25,28,30\}$
(ii) three different values of $p$ depending on the case and (iii) $\mu=0.7$, $\delta=0.9\mu$ and six different $\lambda$-vectors.
These $\lambda$-vectors are the well-known
Median $\lambda=(1,\ldots,1)$,
Anti-$(k_1+k_2)$-trimmed-mean $\lambda=(1,\stackrel{k_1}{\ldots},1,0,\ldots,0,1,\stackrel{k_2}{\ldots},1)$,  $(k_1+k_2)$-Trimmed-mean
$\lambda=(0,\stackrel{k_1}{\ldots},0,1,\ldots,1,0,\stackrel{k_2}{\ldots},0)$,  with $k_1=k_2=\lceil 0.2 n\rceil$,
Center $\lambda=(0,\ldots,0,1)$,
 and $k$-Centrum $\lambda=(0,\ldots,0,1,\stackrel{k}{\ldots},1)$ with $k=\lceil 0.2n \rceil$.
As well as a  $\{0,1\}$-blocks $\lambda$-vector  (three alternate $\{0-1\}$-blocks of lambda weights, i.e.
$\lambda=(0,\ldots,0,1,\ldots,1,0,\ldots,0,\linebreak1,\ldots,1,0,\ldots,0,1,\ldots,1)$).
Therefore, for the each combination of $n$, $p$ and $\lambda$ we have tested five instances. This is, a total number of $450$ problems have been used to test the performance of the proposed models.

\begin{figure*}[t]
\mbox{\hspace*{-3cm}
\includegraphics[scale=0.48]{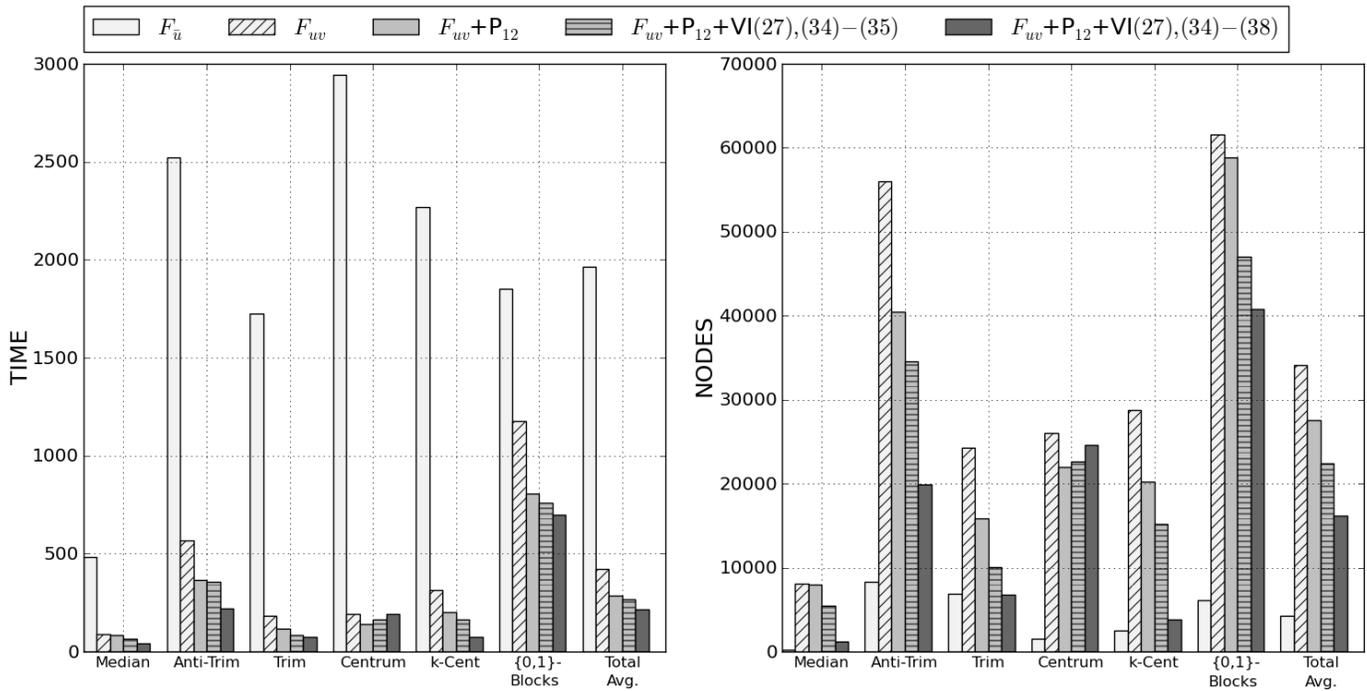}
}
\caption{\footnotesize Summary of computational results carried out in the paper.}
\label{fig1}
\end{figure*}

The capacities were randomly generated in
$\displaystyle [\min_{i} W_i, {(1/2)} \sum_{i \in A} W_i]$.
This generation procedure does not ensure in all cases feasible instances, as
capacity constraints can be very tight for problems with a low number of hubs ($n=3,5$).
Overall, in our experiments we got,  initially 10 infeasible instances out of 75 (13.3\%).
These instances were replaced by new feasible ones (generated with the same capacity
structure). Hence, the reader may observe that the generation procedure gives tight capacity constraints.

First of all, and for the sake of readability, we present in Figure \ref{fig1} a summary of our computational results. A full description of those results is also included in Tables \ref{tabla1}-\ref{tabla6}.

Figure \ref{fig1} shows average results for each one of the considered problem types (different $\lambda$-vectors).
The left chart refers to average CPU times and the right chart to the number of explored nodes of the B\&B tree.
Both charts contain the same number of blocks standing for the different types of $\lambda$-vectors plus and additional block (Total Avg.) for the consolidated average of all $\lambda$-vectors.
Each block compares the behavior of the different formulations ($F_{\bar{u}}$,$F_{uv}$) and their strengthening (variable fixing and valid inequalities).
The heading explains the meaning of bars within each block: $F_{\bar{u}}$ and $F_{uv}$ stand for the corresponding formulations; $F_{uv}$+P$_{12}$ when the two preprocesses $P_{1}(h)$ and $P_{2}(h)$ are applied;  and $F_{uv}$+P$_{12}$+VI when the Valid Inequalities are added as well, denoted with their corresponding references.

Analyzing this figure we observe the improvement obtained with $F_{uv}$ and its strengthening as compared with $F_{\bar{u}}$ or even $F_{uv}$ alone. Actually, the overall reduction in running time with respect to the initial formulation is above $89\%$.

Next, focusing in the best model, namely $F_{uv}+P_{12}+$VI\eqref{in:DV3_3},\eqref{in:DV3c_1}-\eqref{const3}, we observe that the most time consuming problem is the one with $\lambda$-vector given by $\{0,1\}$-blocks with a significant difference with respect to the remaining $\lambda$-vectors. The second most time consuming problem corresponds with the Antitrimmean. Similar conclusions, regarding the number of nodes, are obtained looking at the right chart of
Figure \ref{fig1}.  It is worth mentioning that although the formulation $F_{\bar{u}}$ provides the worst computational times, it is the one that reports the lowest number of nodes in the B\&B tree.  This fact is explained because $F_{\bar u}$ has a larger number of variables than $F_{uv}$ which results in more difficult LP relaxations in each node of the B\&B tree.

{ In spite of that} Figure \ref{fig1} shows the general overview of our computational results, we also include in the following a more detailed analysis based on Tables \ref{tabla1}$-$\ref{tabla6}.

Tables \ref{tabla1} and \ref{tabla2} report the results of the formulations and different preprocessing phases developed in this paper for the  CSA-OMHLP.  The first column of these tables  describes the different $\lambda$-vectors  in the study, the second and third columns report the size of the instances and the
number of hubs to be located, respectively.
The following three columns correspond to some of the computational results obtained by solving the CSA-OMHLP with  $F_{\bar{u}}$ formulation. The next three columns correspond to Formulation  $F_{uv}$, and the rest of the columns to the latter  formulation plus the two preprocessing procedures, i.e $F_{uv}$+P$_{1}$+P$_{2}$.
Columns \textit{RGAP}, \textit{Nodes} and \textit{Time}
stand for the averages of: the gap in the root node, number of nodes in the B\&B tree and
the CPU time in seconds; the time was limited to two hours of CPU.
To obtain a general idea of the comparisons among these averaged values,
for the results in the column {\it Time} and  for different formulations and/or valid inequalities applied, we have accounted the value 7200 seconds for those instances that exceed the time limit. A superindex in their corresponding averaged time value states the number of  instances  exceeding the CPU time limit;  in the same way, the values used to computed the average of the column {\it Nodes} have been the  number of nodes of the B\&B tree when the CPU time limit was reached.

The column \textit{Fixed}  gives the percentage of variables that have been fixed after the  Preprocessing Phases 1 and 2.
Column \textit{Cuts} provides the number of the lower and upper bounds over the {variables $v$}  added to  the model, after running the corresponding preprocessing phases. Finally, \textit{T. prep.} reports  the CPU time in seconds of the corresponding preprocessing phases and column \textit{T. total } reports the overall CPU time in seconds to solve the problem including the corresponding preprocessing phase.
The rows \textit{Average} provide the averaged results among all the tested instances  for each problem and \textit{TOTAL} with respect to the overall set of instances.

Tables  \ref{tabla1} and \ref{tabla2} show, as a general trend, that  the CPU time increases similarly, for all choices of the $\lambda$-vector, with the size of the instances.
We can see that  69 instances required more than two hours to be solved using Formulation
$F_{\bar{u}}$, however  the remaining two analyzed ways  to solve this problem were able to solve all the instances within the CPU  time limit,  except in a few cases for the Antitrimmean and $\{0,1\}$-blocks problem types.
Regarding the running times for the  different types of problems, we can see that Anti-TrimMean and
Blocks  have been the problems that need more time to  be solved.
In any case, we observed that there is a considerable reduction of the running times from the original formulation to the improved formulation, and after applying the  preprocessing phases. From tables \ref{tabla1} and \ref{tabla2},  one can remark that,  Formulation $F_{uv}$ with  Preprocessing Phases 1 and 2 provides better results, with a reduction of  around a 86\% of the time with respect  to Formulation $F_{\bar{u}}$  (taking into account that the latter formulation was not able to solve all the studied instances before the time limit was exceeded). Moreover, the reduction of the running times with respect to Formulation
$F_{uv}$ (without preprocessing phases) is around 32 \%.

As for the comparison between  RGAP's, we observe that the average gap of the linear relaxation after  preprocessing reduces around  32\% from the original formulation
(Formulation $F_{\bar{u}}$) to Formulation
$F_{uv}$ with Preprocessing Phases 1 and 2.
In any case, it is also worth
noting that the gap from the original formulation to the improved formulation (without any
preprocessing phase) is reduced around  9\%, what implies that even though this improved formulation uses much less number of variables, it provides a better RGAP.

Tables \ref{table3} and \ref{table4} present several improvements to Formulation  $F_{uv}$ with Preprocessing Phases 1 and 2.  In particular, the second, third and fourth blocks of columns summarize the results when  the combination of  valid inequalities
(\ref{in:DV3_3})-(\ref{in:DV3_2});
(\ref{in:DV3_3}),(\ref{in:DV3b_1})-(\ref{in:DV3b_2}) and
(\ref{in:DV3_3}),(\ref{in:DV3c_1})-(\ref{in:DV3c_2})
are added to this formulation, respectively. In general,  we can see that the latter combination provides the best results. In particular, there is an improvement of 7 \% of the running
times with respect to Formulation { $F_{uv}$} with  Preprocessing Phases 1 and 2. Moreover, with respect to the number of nodes the improvement is around 15\%, but the RGAP is similar for all the analyzed reinforcements.

Since the best behavior observed was obtained with reinforcement given by (\ref{in:DV3_3}), (\ref{in:DV3c_1})-(\ref{in:DV3c_2}), in the rest of our tests we have used this configuration to make further strengthening.
Tables \ref{table5} and \ref{table6} present several improvements to the
Formulation { $F_{uv}$} with  Preprocessing Phases 1 and 2 and valid
inequalities (\ref{in:DV3_3}), (\ref{in:DV3c_1})-(\ref{in:DV3c_2}). In particular, the second, third and fourth blocks of columns report the results after including the family of valid inequalities
(\ref{in:cortepacking1}), (\ref{in:cortepacking1})-(\ref{const1}) and
(\ref{in:cortepacking1})-(\ref{const3}), respectively.
These tables show that the best results are obtained when all  valid inequalities,
(\ref{in:cortepacking1})-(\ref{const3}), are added to the current configuration; with an improvement of 18\% in the running time, 28\% in the
number of nodes and 12 \% in the RGAP.

The overall conclusion of our experiment is that in order to solve CSA-OMHLP the best combination of formulation and strengthening is to use $F_{uv}+P_{12}+$(\ref{in:DV3_3})$+$(\ref{in:DV3c_1})$+$(\ref{in:DV3c_2})$+$(\ref{in:cortepacking1})$+$(\ref{const3}). This configuration allows to solve medium size instances within 10 minutes of CPU time.

\section{Concluding remarks}\label{Sec:Concl}

This paper can be  considered as an initial attempt  to address the  capacitated single-allocation
ordered hub location problem. The formulations, strengthening and preprocessing phases developed in this paper  provide a promising approach to solve  the above mentioned  problem
although { so far} only medium size problems are reasonably well-solved. Thus, this work opens interesting possibilities  to study-and-develop ad-hoc solution procedures that allow us to consider larger size instances of this problem.  Moreover, it also  points out the possibility of developing heuristics approaches that will give good solutions in competitive running times.
All in all, this paper   shows the usefulness of using  covering formulations and their corresponding strengthening for solving capacitated versions of ordered hub location problems.

\section*{Acknowledgement}

This research has been partially supported by Spanish Ministry of Education and Science/FEDER grants numbers MTM2010-19576-C02-(01-02), MTM2013-46962-C02-(01-02),  Junta de Andaluc{\'\i}a grant number FQM 05849 and Fundaci\'on S\'eneca, grant number 08716/PI/08.


\addcontentsline{toc}{chapter}{Bibliography}
\bibliographystyle{amsalpha}

\providecommand{\bysame}{\leavevmode\hbox to3em{\hrulefill}\thinspace}
\providecommand{\MR}{\relax\ifhmode\unskip\space\fi MR }
\providecommand{\MRhref}[2]{%
  \href{http://www.ams.org/mathscinet-getitem?mr=#1}{#2}
} \providecommand{\href}[2]{#2}

\begin{table}[htbp]
\vspace*{-1.cm}
\caption{Computational results for Formulations   $F_{\bar{u}}$
vs $F_{uv}$  with preprocessing.
\label{tabla1} }
\centering
{
\tiny \setlength{\tabcolsep}{0.2cm}
\hspace*{-0.5cm}
\begin{tabular}{|c|c|c|rrr|rrr|rrrrrr|}
\hline
           & \multirow{ 2}{*}{n} & \multirow{ 2}{*}{p} & \multicolumn{ 3}{|c|}{{\bf Formulation $F_{\bar{u}}$ 
           }} &      \multicolumn{ 3}{|c|}{{\bf Formulation   $F_{uv}$ 
           }} &                                          
          \multicolumn{ 6}{|c|}{{\bf Formulation   $F_{uv}$ +
{ P$_{1}$+P$_{2}$}
}} \\
\cline{4-15}
           & &  &       RGAP &      Nodes &       Time &       RGAP &      Nodes &       Time &       RGAP &      Nodes &    T total &    T. prep &      Fixed &       Cuts  \\
\hline
\multirow{15}{*}{\begin{sideways}\texttt{MEDIAN}\end{sideways}}  &      15 &          3 &      17,56 &       17 &      7,9 &      17,56 &      240 &      1,8 &      15,16 &      195 &        1,4 &        1,4 &        2,4 &        240 \\

\multicolumn{ 1}{|c|}{} &         15 &          5 &      7,60 &        3 &      3,8 &      7,60 &      230 &     0,8 &      6,78 &      177 &        0,7 &        1,2 &        1,4 &      180 \\

\multicolumn{ 1}{|c|}{} &         15 &          8 &      8,49 &        1 &      1,0 &      8,49 &      259 &      0,6 &      4,39 &      185 &        0,4 &        1,0 &        0,8 &      133 \\

\multicolumn{ 1}{|c|}{} &         20 &          3 &      19,54 &       75 &      69,5 &      19,54 &      505 &      6,3 &      16,31 &      472 &        5,8 &        3,1 &        3,9 &      405 \\

\multicolumn{ 1}{|c|}{} &         20 &          8 &      6,05 &        1 &      3,9 &      6,05 &     1024 &      4,2 &      5,49 &     1122 &        3,8 &        2,2 &        1,9 &      215 \\

\multicolumn{ 1}{|c|}{} &         20 &         10 &      2,88 &        1 &      3,1 &      2,88 &      591 &      1,9 &      2,59 &      560 &        1,6 &        1,7 &        1,8 &        181 \\

\multicolumn{ 1}{|c|}{} &         25 &          3 &      19,00 &    1.051 &      509,5 &         19,00 &      478 &      21,2 &      15,21 &      408 &       18,5 &        8,9 &        7,2 &      800 \\

\multicolumn{ 1}{|c|}{} &         25 &          8 &      6,77 &        1 &      12,7 &      6,77 &     4239 &      30,9 &      5,95 &     6177 &       33,6 &        5,2 &        4,1 &      459 \\

\multicolumn{ 1}{|c|}{} &         25 &         10 &       13,10 &       15 &      48,6 &       13,10 &     6656 &      33,4 &      6,31 &     4686 &       26,0 &        5,2 &        3,4 &      435 \\

\multicolumn{ 1}{|c|}{} &         28 &          3 &      21,86 &      270 &      322,7 &      21,86 &     2052 &      66,8 &      16,01 &     2246 &       76,6 &       13,2 &       10,5 &      960 \\

\multicolumn{ 1}{|c|}{} &         28 &          8 &      9,21 &      143 &      415,6 &      9,21 &    16820 &      165,5 &      8,15 &    10560 &       97,4 &        7,7 &        6,5 &      597 \\

\multicolumn{ 1}{|c|}{} &         28 &         10 &       7,36 &      129 &      321,3 &       7,36 &    27440 &      193,1 &      6,78 &    32680 &      243,1 &        6,7 &        3,8 &      396 \\

\multicolumn{ 1}{|c|}{} &         30 &          3 &      32,36 &    1.416 &    4007,0$^{(1)}$ &      32,32 &     3701 &      167,7 &      20,92 &     2564 &      127,8 &       18,8 &       11,0 &       1128 \\

\multicolumn{ 1}{|c|}{} &         30 &          8 &      18,69 &      360 &       1346,0 &      18,69 &    39430 &      455,5 &      10,86 &    36790 &      437,3 &       12,8&        6,2 &      817 \\

\multicolumn{ 1}{|c|}{} &         30 &         10 &       9,18 &       22 &      186,3 &       9,18 &    18420 &      176,8 &       6,98 &    20130 &      207,6 &       11,4 &        5,4 &      685 \\
\cline{2-15}\rowcolor[gray]{0.9}
\multicolumn{ 1}{|c|}{} & \multicolumn{ 2}{c|}{{\bf Average}} &{\bf 13,30} & {\bf 234} & {\bf 483,9} & {\bf 13,30} & {\bf 8139} & {\bf 88,4} &  {\bf 9,90} & {\bf 7930} & {\bf 85,4} &  {\bf 6,7} &  {\bf 4,7} & {\bf 509} \\
\hline
\multirow{15}{*}{\begin{sideways}\texttt{ANTITRIMMEAN}\end{sideways}} &       15 &          3 &      23,75 &      150 &      33,6 &      23,75 &      614 &      2,8 &      17,62 &      645 &        2,0 &        1,4 &        6,1 &      164 \\

\multicolumn{ 1}{|c|}{} &         15 &          5 &      13,93 &       99 &      12,9 &      13,93 &     1039 &      2,2 &      8,38 &      803 &        1,4 &        1,1 &        4,3 &      104 \\

\multicolumn{ 1}{|c|}{} &         15 &          8 &       13,20 &       29 &      9,4 &       13,20 &      618 &      1,1 &      6,82 &      724 &        0,9 &        1,0 &        2,6 &       91 \\

\multicolumn{ 1}{|c|}{} &         20 &          3 &       23,50 &    2448 &        272,0 &       23,50 &     2032 &      16,8 &      17,44 &     1150 &        9,2 &        3,1 &        9,8 &      273 \\

\multicolumn{ 1}{|c|}{} &         20 &          8 &      9,45 &      379 &      79,9 &      9,45 &     4637 &      11,9 &      7,66 &     3009 &        7,8 &        2,2 &        5,2 &      126 \\

\multicolumn{ 1}{|c|}{} &         20 &         10 &      4,66 &       73 &      36,4 &      4,66 &     3150 &      6,9 &      4,28 &     3394 &        6,4 &        1,7 &        4,4 &      123 \\

\multicolumn{ 1}{|c|}{} &         25 &          3 &      24,64 &   12130 &       2037,0 &      24,64 &     1984 &      53,8 &      17,29 &     1218 &       29,0 &        8,9 &       18,9 &        540 \\

\multicolumn{ 1}{|c|}{} &         25 &          8 &      10,32 &      590 &      344,4 &      10,32 &     6411 &      40,2 &      7,91 &     6841 &       33,3 &        5,2 &       10,4 &        324 \\

\multicolumn{ 1}{|c|}{} &         25 &         10 &      17,25 &   13650 &       1217,0 &      17,25 &    27180 &      115,1 &       9,81 &    25390 &      102,1 &        5,2 &        8,6 &        330 \\

\multicolumn{ 1}{|c|}{} &         28 &          3 &      25,94 &   16460 &    5111,0$^{(1)}$ &      25,94 &    26600 &      957,1 &      17,36 &     7432 &      168,4 &       13,2 &       22,1 &      717 \\

\multicolumn{ 1}{|c|}{} &         28 &          8 &      13,58 &   12390 &    5688,0$^{(3)}$ &      13,52 &    42860 &      381,4 &      10,67 &    34800 &      271,1 &        7,8 &       16,8 &      301 \\

\multicolumn{ 1}{|c|}{} &         28 &         10 &      10,95 &   13830 &    3524,0$^{(2)}$ &      10,92 &   262500 &       1595,0 &      9,14 &   131900 &      832,1 &        6,7 &        9,7 &      242 \\

\multicolumn{ 1}{|c|}{} &         30 &          3 &      36,85 &   11110 &    6584,0$^{(4)}$ &      35,55 &    15670 &        591,0 &      21,74 &     7220 &      271,6 &       18,8 &       27,2 &        770 \\

\multicolumn{ 1}{|c|}{} &         30 &          8 &       23,80 &   19960 &    7200,0$^{(5)}$ &      23,72 &   299700 &       3402,0 &      14,48 &   273600 &     2896,0 &       12,8 &       15,8 &      598 \\

\multicolumn{ 1}{|c|}{} &         30 &         10 &      13,53 &   21200 &    5707,0$^{(3)}$ &       13,50 &   145000 &       1310,0 &      10,56 &   109700 &      863,2 &       11,4 &       15,4 &      464 \\
\cline{2-15}\rowcolor[gray]{0.9}
\multicolumn{ 1}{|c|}{} & \multicolumn{ 2}{c|}{{\bf Average}} &{\bf 17,70} & {\bf 8300} & {\bf 2523,8} & {\bf 17,60} & {\bf 56000} & {\bf 565,8} & {\bf 12,10} & {\bf 40522} & {\bf 366,3} &  {\bf 6,7} & {\bf 11,8} & {\bf 345} \\
\hline
\multirow{15}{*}{\begin{sideways}\texttt{TRIMMEAN}\end{sideways}}&     15 &          3 &      19,46 &      177 &       22,9 &      17,43 &     1066 &      2,1 &      16,44 &      729 &        1,4 &        1,4 &        0,8 &       76 \\

\multicolumn{ 1}{|c|}{} &         15 &          5 &      11,26 &       38 &      13,6 &      8,65 &     1209 &      1,1 &      8,58 &      646 &        0,7 &        1,2 &        0,8 &         67 \\

\multicolumn{ 1}{|c|}{} &         15 &          8 &       8,72 &        4 &      3,6 &      7,88 &      510 &      0,4 &      7,88 &      322 &        0,3 &        1,0 &        0,3 &       20 \\

\multicolumn{ 1}{|c|}{} &         20 &          3 &      20,77 &    5452 &      335,9 &      19,38 &     1047 &      7,0 &      17,65 &     1022 &        6,2 &        3,1 &        1,9 &        132 \\

\multicolumn{ 1}{|c|}{} &         20 &          8 &       12,10 &    1107 &      96,7 &      9,48 &     6847 &      13,4 &      9,08 &     4512 &        7,9 &        2,2 &        1,2 &       88 \\

\multicolumn{ 1}{|c|}{} &         20 &         10 &      8,88 &       63 &       36,7 &      6,17 &     2309 &      4,4 &      5,97 &     1620 &        2,7 &        1,7 &        0,9 &       58 \\

\multicolumn{ 1}{|c|}{} &         25 &          3 &      19,74 &   18180 &       2157,0 &      18,78 &     2039 &      31,4 &      16,87 &     2081 &       25,1 &        8,9 &        3,2 &      260 \\

\multicolumn{ 1}{|c|}{} &         25 &          8 &      10,55 &    5365 &        804,0 &      8,94 &    61330 &      211,7 &      8,63 &    21460 &       74,2 &        5,2 &        2,2 &      140 \\

\multicolumn{ 1}{|c|}{} &         25 &         10 &      13,58 &    1186 &      225,4 &      12,34 &     2642 &      9,6 &      12,28 &     1760 &        6,5 &        5,2 &        1,9 &        117 \\

\multicolumn{ 1}{|c|}{} &         28 &          3 &       21,50 &   18050 &    4959,0$^{(1)}$ &      20,54 &     8611 &      149,2 &      17,49 &     3800 &       69,6 &       13,3 &        3,9 &        243 \\

\multicolumn{ 1}{|c|}{} &         28 &          8 &      12,23 &   12740 &       2304,0 &      10,74 &    47430 &      309,6 &      10,24 &    20290 &      130,9 &        7,7 &        3,8 &      234 \\

\multicolumn{ 1}{|c|}{} &         28 &         10 &      11,27 &    6731 &       1467,0 &      9,25 &    35750 &      178,7 &      9,171 &    14770 &       69,3 &        6,7 &        2,8 &      174 \\

\multicolumn{ 1}{|c|}{} &         30 &          3 &      26,67 &    9820 &    6506,0$^{(4)}$ &      25,15 &     9904 &      294,1 &      22,08 &     3302 &      124,4 &       18,8 &        5,6 &      358 \\

\multicolumn{ 1}{|c|}{} &         30 &          8 &      17,04 &   10390 &    3766,0$^{(1)}$ &      15,44 &    90810 &      868,1 &      15,09 &    54920 &      483,1 &       12,8 &        3,6 &      243 \\

\multicolumn{ 1}{|c|}{} &         30 &         10 &      11,37 &   14330 &       3181,0 &      9,44 &    92100 &        690,0 &       9,40 &   107100 &      757,0 &       11,4 &        3,2 &      227 \\
\cline{2-15}\rowcolor[gray]{0.9}
\multicolumn{ 1}{|c|}{} & \multicolumn{ 2}{c|}{{\bf Average}} &   {\bf 15,00} & {\bf 6909} & {\bf 1725,3} & {\bf 13,30} & {\bf 24240} & {\bf 184,7} & {\bf 12,50} & {\bf 15889} & {\bf 117,3} &  {\bf 6,7} &  {\bf 2,4} & {\bf 163} \\
\hline
\end{tabular} }
\end{table}

\begin{table}[htbp]
\vspace*{-0.5cm}
\caption{Computational results for Formulations  $F_{\bar{u}}$
vs $F_{uv}$ with preprocessing(II).
\label{tabla2} }
\centering
{
\tiny \setlength{\tabcolsep}{0.2cm}
\hspace*{-1cm}
\begin{tabular}{|c|c|c|rrr|rrr|rrrrrr|}
\hline
           & \multirow{ 2}{*}{n} & \multirow{ 2}{*}{p} & \multicolumn{ 3}{|c|}{{\bf Formulation $F_{\bar{u}}$ 
           }} &      \multicolumn{ 3}{|c|}{{\bf Formulation   $F_{uv}$ 
           }} &                                           
          \multicolumn{ 6}{|c|}{{\bf Formulation   $F_{uv}$ +
{ P$_{1}$+P$_{2}$}
}} \\
\cline{4-15}
           & &  &       RGAP &      Nodes &       Time &       RGAP &      Nodes &       Time &       RGAP &      Nodes &    T total &    T. prep &      Fixed &       Cuts  \\
\hline
\multirow{15}{*}{\begin{sideways}\texttt{CENTER}\end{sideways}} &    15 &          3 &      23,84 &      159&       36,7 &      22,87 &      780 &      2,5 &       16,90 &      489 &     0,9 &      1,4 &      5,1 &          0 \\

\multicolumn{ 1}{|c|}{} &         15 &          5 &      16,54 &      194 &      24,5 &      14,97 &     1374 &       1,9 &      8,48 &      808 &     0,9 &      1,1 &      4,7 &          0 \\

\multicolumn{ 1}{|c|}{} &         15 &          8 &      19,85 &       84 &       21,1 &      16,45 &      622 &      1,0 &      8,83 &      501,0 &     0,6 &      1,0 &      3,9 &          0 \\

\multicolumn{ 1}{|c|}{} &         20 &          3 &      21,36 &    1.225 &      347,5 &      20,76 &     1118 &      9,9 &      16,72 &    1017 &      4,8 &      3,2 &      8,9 &          0 \\

\multicolumn{ 1}{|c|}{} &         20 &          8 &       14,30 &      775 &      224,6 &      12,48 &     4195 &      8,9 &      11,58 &    3712 &      6,2 &      2,2 &      4,9 &          0 \\

\multicolumn{ 1}{|c|}{} &         20 &         10 &      12,16 &      438 &      88,6 &      9,32 &     2675 &      5,9 &      8,817 &    3019 &      5,3 &      1,7 &       3,8 &          0 \\

\multicolumn{ 1}{|c|}{} &         25 &          3 &      21,97 &    1785 &       1626,0 &      21,58 &     1948 &      31,2 &      14,57 &      938 &      12,3 &      8,9 &      18,4 &          0 \\

\multicolumn{ 1}{|c|}{} &         25 &          8 &       15,70 &    2763 &       1663,0 &       14,50 &    45250 &      161,6 &      11,33 &   26520 &      87,4 &      5,2 &      10,7 &          0 \\

\multicolumn{ 1}{|c|}{} &         25 &         10 &      27,76 &    1901 &       1525,0 &      25,91 &     5302 &      21,2 &      16,16 &    5134 &      17,7 &      5,2 &      11,1 &          0 \\

\multicolumn{ 1}{|c|}{} &         28 &          3 &       22,50 &    2698 &    5939,0$^{(2)}$ &       21,30 &     5756 &        121,0 &      17,56 &    4172 &      66,8 &      13,2 &       18,5 &          0 \\

\multicolumn{ 1}{|c|}{} &         28 &          8 &      18,35 &    1931 &    6031,0$^{(2)}$ &      16,63 &    30380 &        198,0 &      15,77 &   29450 &      170,1 &      7,8 &      11,5 &          0 \\

\multicolumn{ 1}{|c|}{} &         28 &         10 &      16,65 &    3266 &    6621,0$^{(3)}$ &      15,09 &    54930 &      272,6 &      14,53 &   48120 &      218,4 &      6,7 &      8,6 &          0 \\

\multicolumn{ 1}{|c|}{} &         30 &          3 &       42,50 &    1725 &    6234,0$^{(4)}$ &      31,02 &     8074 &      293,5 &      19,75 &    6536 &      147,3 &       18,7 &      24,4 &          0 \\

\multicolumn{ 1}{|c|}{} &         30 &          8 &      32,41 &    1035 &    7200,0$^{(5)}$ &      27,59 &   134300 &       1096,0 &      14,63 &  103700 &      786,6 &      12,8 &      21,2 &          0 \\

\multicolumn{ 1}{|c|}{} &         30 &         10 &      17,93 &    2533 &    6563,0$^{(4)}$ &       16,20 &    93330 &      646,2 &      11,77 &   96140 &      593,6 &      11,4 &      17,3 &          0 \\
\cline{2-15}\rowcolor[gray]{0.9}
\multicolumn{ 1}{|c|}{} & \multicolumn{ 2}{c|}{{\bf Average}} & {\bf 21,60} & {\bf 1501} & {\bf 2945,8} & {\bf 19,10} & {\bf 26002} & {\bf 191,4} & {\bf 13,80} & {\bf 22017} & {\bf 141,3} &  {\bf 6,7} & {\bf 11,5} &  {\bf 0} \\
\hline
\multirow{15}{*}{\begin{sideways}\texttt{K-CENTRUM}\end{sideways}}&     15 &          3 &      26,74 &      148 &      39,9 &      23,07 &     1031 &      3,0 &       17,60 &      588 &       1,8 &       1,4 &      1,5 &      143 \\

\multicolumn{ 1}{|c|}{} &         15 &          5 &      16,81 &        7 &      14,5 &      10,97 &     1093 &      1,9 &      7,91 &      498 &      1,0 &      1,1 &      1,3 &      161 \\

\multicolumn{ 1}{|c|}{} &         15 &          8 &      15,68 &        1 &      3,1 &      8,49 &      259 &      0,6 &      4,39 &      185 &     0,4 &      1,0 &     0,8 &      133 \\

\multicolumn{ 1}{|c|}{} &         20 &          3 &      27,98 &    1995 &      526,5 &      25,45 &     2426 &         17,0 &      18,23 &    1528 &      9,1 &      3,1 &       2,6 &      247 \\

\multicolumn{ 1}{|c|}{} &         20 &          8 &      14,46 &      841 &      146,9 &      9,49 &     5136 &      12,7 &      7,12 &    2902 &      6,2 &      2,2 &      1,7 &      144 \\

\multicolumn{ 1}{|c|}{} &         20 &         10 &      9,53 &      356 &      67,7 &      4,84 &     4255 &      7,7 &      4,04 &    3705 &      5,9 &      1,7 &       1,2 &      127 \\

\multicolumn{ 1}{|c|}{} &         25 &          3 &      29,31 &    4964 &       2562,0 &      27,24 &     2727 &      49,1 &      17,48 &    2693 &      32,3 &       8,9 &      5,3 &        442 \\

\multicolumn{ 1}{|c|}{} &         25 &          8 &      14,14 &       52 &      476,7 &      9,79 &     5928 &       38,5 &      7,33 &    3910 &      23,2 &      5,2 &      2,7 &      342 \\

\multicolumn{ 1}{|c|}{} &         25 &         10 &      20,63 &    1976 &      512,9 &      14,78 &    13330 &      64,9 &      7,34 &   12820 &      52,6 &      5,2 &        2,2 &        358 \\

\multicolumn{ 1}{|c|}{} &         28 &          3 &      29,78 &    2991 &    6343,0$^{(3)}$ &      27,99 &    13280 &      264,2 &      17,17 &    4971 &      111,3 &      13,3 &       6,2 &        466 \\

\multicolumn{ 1}{|c|}{} &         28 &          8 &      15,97 &    9132 &    4372,0$^{(1)}$ &      12,12 &    60840 &      531,5 &      9,44 &   21180 &      178,9 &      7,8 &      3,7 &      431 \\

\multicolumn{ 1}{|c|}{} &         28 &         10 &      14,09 &    6460 &    2576,0$^{(1)}$ &      9,11 &    73800 &      526,4 &      7,63 &   35460 &      228,2 &      6,7 &       2,7 &      324 \\

\multicolumn{ 1}{|c|}{} &         30 &          3 &      41,03 &      671 &    7200,0$^{(5)}$ &      36,42 &    18830 &      616,8 &      21,74 &    4588 &      187,6 &      18,7 &      7,3 &      720 \\

\multicolumn{ 1}{|c|}{} &         30 &          8 &      25,04 &    5859 &    6632,0$^{(4)}$ &      21,06 &   174500 &       2047,0 &       12,20 &  171500 &       1830,0 &      12,8 &      4,4 &        695 \\

\multicolumn{ 1}{|c|}{} &         30 &         10 &      14,74 &    2560 &       2536,0 &      10,24 &    53280 &      506,9 &      7,63 &   37330 &        340,0 &      11,4 &      4,8 &      648 \\
\cline{2-15}\rowcolor[gray]{0.9}
\multicolumn{ 1}{|c|}{} & \multicolumn{ 2}{c|}{{\bf Average}} & {\bf 21,10} & {\bf 2534} & {\bf 2267,7} & {\bf 16,70} & {\bf 28714} & {\bf 312,5} & {\bf 11,20} & {\bf 20257} & {\bf 200,6} &  {\bf 6,7} &  {\bf 3,2} & {\bf 359} \\
\hline
\multirow{15}{*}{\begin{sideways}\texttt{BLOCKS}\end{sideways}} &            15 &          3 &      19,46 &       94 &      15,5 &      19,46 &     2210 &      8,9 &      15,97 &    1729 &      3,8 &      1,45 &      8,7 &       91 \\

\multicolumn{ 1}{|c|}{} &         15 &          5 &      12,24 &       33 &      10,9 &      11,79 &     1122 &       3,7 &      8,08 &      860 &       2,4 &      1,2 &      6,7 &      117 \\

\multicolumn{ 1}{|c|}{} &         15 &          8 &      12,71 &        9 &      7,6 &      12,37 &     1065 &      1,8 &       6,25 &      715 &      1,1 &      1,0 &      3,3 &         76 \\

\multicolumn{ 1}{|c|}{} &         20 &          3 &      21,09 &      606 &        113,0 &      20,86 &     3721 &      30,7 &      17,58 &    5777 &       24,7 &      3,2 &       16,1 &      219 \\

\multicolumn{ 1}{|c|}{} &         20 &          8 &      8,55 &      402 &      67,3 &      8,55 &    12670 &      41,3 &      7,46 &    8064 &      20,4 &      2,2 &      4,8 &      116 \\

\multicolumn{ 1}{|c|}{} &         20 &         10 &      4,70 &       97 &       33,2 &      4,32 &     3590 &      8,4 &      4,00 &    4585 &      9,0 &      1,7 &      3,7 &      122 \\

\multicolumn{ 1}{|c|}{} &         25 &          3 &      21,51 &    4850 &      983,7 &      21,27 &     5637 &      114,4 &      16,45 &    4087 &      59,1 &      8,9 &      27,5 &      484 \\

\multicolumn{ 1}{|c|}{} &         25 &          8 &      9,99 &    1675 &      429,2 &      9,94 &    13400 &        115,0 &      8,12 &   12050 &      89,3 &      5,2 &      17,4 &        260 \\

\multicolumn{ 1}{|c|}{} &         25 &         10 &      17,57 &    6874 &       1021,0 &       17,40 &    46270 &      216,9 &      10,01 &   40870 &      177,5 &      5,2 &      7,4 &        268 \\

\multicolumn{ 1}{|c|}{} &         28 &          3 &      22,45 &    5627 &       2887,0 &      22,45 &     9111 &      472,4 &      17,16 &    8569 &      267,1 &      13,4 &      48,3 &      482 \\

\multicolumn{ 1}{|c|}{} &         28 &          8 &      10,44 &   14490 &       2687,0 &      10,25 &    96840 &       1515,0 &       9,32 &   54910 &        745,0 &      7,7 &      24,1 &      357 \\

\multicolumn{ 1}{|c|}{} &         28 &         10 &      9,39 &   13040 &    3181,0$^{(1)}$ &      9,26 &   233200 &       2683,0 &      8,67 &  184700 &       1601,0 &      6,7 &      15,4 &      218 \\

\multicolumn{ 1}{|c|}{} &         30 &          3 &      33,18 &   10180 &    5755,0$^{(2)}$ &      33,04 &    10250 &      753,4 &      22,39 &   12840 &      703,5 &      18,9 &      56,4 &      664 \\

\multicolumn{ 1}{|c|}{} &         30 &          8 &      21,41 &   15590 &    6438,0$^{(4)}$ &      21,46 &   201400 &       6587,0 &      12,95 &  193300 &    5016,0$^{(3)}$ &      12,9 &      24,5 &        511 \\

\multicolumn{ 1}{|c|}{} &         30 &         10 &       11,90 &   18420 &    4144,0$^{(2)}$ &      11,66 &   283700 &       5113,0 &      9,45 &  256600 &    3402,0$^{(1)}$ &      11,51 &      23,28 &      465 \\
\hline
\cline{2-15}\rowcolor[gray]{0.9}
\multicolumn{ 1}{|c|}{} & \multicolumn{ 2}{c|}{{\bf Average}} & {\bf 15,80} & {\bf 6132} & {\bf 1851,6} & {\bf 15,60} & {\bf 61612} & {\bf 1177,7} & {\bf 11,60} & {\bf 52644} & {\bf 808,1} &  {\bf 6,8} & {\bf 19,2} & {\bf 297} \\
\hline \rowcolor[gray]{0.8}
    {\bf } & \multicolumn{ 2}{|c|}{{\bf TOTAL}} & {\bf 17,40} & {\bf 4268} & {\bf 1966,3} & {\bf 15,90} & {\bf 34118} & {\bf 420,1} & {\bf 11,80} & {\bf 26543} & {\bf 286,5} &  {\bf 6,7} &  {\bf 8,8} & {\bf 279} \\
\hline
\end{tabular}}
\end{table}


\begin{table}[htbp]
\vspace*{-1.cm}
 \caption{Computational results for Formulation$F_{uv}$ with Valid Inequalities (I) \label{table3}}
 \centering
{
\tiny \setlength{\tabcolsep}{0.2cm}
\vspace*{-0.3cm}
\hspace*{-1cm}

\begin{tabular}{|c|c|c|rrr|rrr|rrr|rrr|}
\hline
  &      \multirow{ 2}{*}{n} & \multirow{ 2}{*}{p}&
\multicolumn{ 3}{c|}{{ P$_{1}$+P$_{2}=$  P$_{12}$ }} &
\multicolumn{ 3}{c|}{{\bf P$_{12}$+V.I.(\ref{in:DV3_3})-(\ref{in:DV3_2})}} &
\multicolumn{ 3}{c|}{{\bf P$_{12}$+V.I.(\ref{in:DV3_3}),(\ref{in:DV3b_1})-(\ref{in:DV3b_2})}} &
\multicolumn{ 3}{c|}{{\bf P$_{12}$+V.I.(\ref{in:DV3_3}),(\ref{in:DV3c_1})-(\ref{in:DV3c_2})}} \\
\cline{4-15}
           &           &          &       RGAP &      Nodes &       time &       RGAP &      Nodes &       time &       RGAP &      Nodes &       time &       RGAP &      Nodes &       time \\
 \hline
\multirow{15}{*}{\begin{sideways}\texttt{MEDIAN}\end{sideways}} &         15 &          3 &      15,16 &      195 &      1,4 &      14,21 &      199 &      1,4 &      14,21 &      262 &      1,5 &      14,21 &      239 &      1,4 \\

\multicolumn{ 1}{|c|}{} &         15 &          5 &      6,78 &      177 &     0,7 &        6,50 &      154 &     0,7 &        6,50 &      232 &       0,8 &        6,50 &      230 &     0,8 \\

\multicolumn{ 1}{|c|}{} &         15 &          8 &      4,39 &      185 &     0,4 &      4,33 &      110 &     0,4 &      4,33 &      157 &     0,4 &      4,34 &      162 &     0,5 \\

\multicolumn{ 1}{|c|}{} &         20 &          3 &      16,31 &      472 &       5,8 &      15,66 &      379 &      7,1 &      15,66 &      540 &      7,4 &      15,66 &      560 &      7,3 \\

\multicolumn{ 1}{|c|}{} &         20 &          8 &      5,49 &    1122 &      3,8 &      4,98 &      981 &      4,1 &      4,98 &    1189 &      4,4 &      4,98 &      936 &      3,9 \\

\multicolumn{ 1}{|c|}{} &         20 &         10 &      2,59 &      560 &      1,6 &      2,33 &      658 &      1,9 &      2,33 &      654 &      1,8 &      2,33 &      857 &      2,4 \\

\multicolumn{ 1}{|c|}{} &         25 &          3 &      15,21 &      408 &      18,5 &      14,79 &      502 &      21,3 &      14,79 &      490 &      20,9 &      14,79 &      390 &      19,9 \\

\multicolumn{ 1}{|c|}{} &         25 &          8 &      5,95 &    6177 &       33,6 &      5,39 &     1772 &      16,8 &      5,39 &    2994 &      21,6 &      5,39 &     3301 &      26,5 \\

\multicolumn{ 1}{|c|}{} &         25 &         10 &      6,31 &    4686 &         26,0 &       5,97 &     2748 &      16,1 &       5,97 &    3141 &       18,5 &       5,97 &     2734 &       16,8 \\

\multicolumn{ 1}{|c|}{} &         28 &          3 &      16,01 &    2246 &      76,6 &      15,72 &     1968 &      72,6 &      15,72 &    1871 &      71,5 &      15,72 &     2055 &      71,9 \\

\multicolumn{ 1}{|c|}{} &         28 &          8 &      8,15 &   10560 &      97,4 &      7,53 &     5937 &      62,5 &      7,53 &    5999 &      66,2 &      7,53 &     6398 &      70,4 \\

\multicolumn{ 1}{|c|}{} &         28 &         10 &      6,78 &   32680 &      243,1 &      6,48 &    21080 &      155,2 &      6,48 &   16800 &      136,4 &      6,48 &    17930 &      125,4 \\

\multicolumn{ 1}{|c|}{} &         30 &          3 &      20,92 &    2564 &      127,8 &       20,60 &     2300 &      132,9 &       20,60 &    2797 &      149,9 &       20,60 &     2770 &      142,3 \\

\multicolumn{ 1}{|c|}{} &         30 &          8 &      10,86 &   36790 &      437,3 &      10,35 &    31990 &        448,0 &      10,35 &   25920 &      338,3 &      10,35 &    25060 &      332,7 \\

\multicolumn{ 1}{|c|}{} &         30 &         10 &       6,98 &   20130 &      207,6 &      6,75 &    17230 &      167,8 &      6,75 &   15720 &      160,4 &      6,75 &    18130 &      179,8 \\
\cline{2-15}\rowcolor[gray]{0.9}
\multicolumn{ 1}{|c|}{} & \multicolumn{ 2}{c|}{{\bf Average}} &  {\bf 9,90} & {\bf 7930} & {\bf 85,4} &  {\bf 9,40} & {\bf 5867} & {\bf 73,9} &  {\bf 9,40} & {\bf 5251} & {\bf 66,7} &  {\bf 9,40} & {\bf 5450} & {\bf 66,8} \\
\hline
\multirow{15}{*}{\begin{sideways}\texttt{ANTI-TRIMMEAN}\end{sideways}} &         15 &          3 &      17,62 &      645 &      2,0 &      17,18 &      559 &      1,9 &      17,18 &      709 &      2,1 &      17,18 &      556 &      2,0 \\

\multicolumn{ 1}{|c|}{} &         15 &          5 &      8,38 &      803 &      1,4 &      8,36 &      823 &      1,6 &      8,36 &      615 &      1,4 &      8,36 &      705 &      1,5 \\

\multicolumn{ 1}{|c|}{} &         15 &          8 &      6,82 &      724 &      0,9 &       6,82 &      808 &     0,9 &       6,82 &      723 &      0,9 &       6,82 &      736 &     1,0 \\

\multicolumn{ 1}{|c|}{} &         20 &          3 &      17,44 &    1150 &       9,2 &         17,00 &     1471 &      11,4 &         17,00 &    1455 &      11,7 &         17,00 &     1370 &      10,9 \\

\multicolumn{ 1}{|c|}{} &         20 &          8 &      7,66 &    3009 &      7,8 &      7,54 &     2684 &      7,1 &      7,54 &    2716 &      7,6 &      7,54 &     2709 &      7,7 \\

\multicolumn{ 1}{|c|}{} &         20 &         10 &      4,28 &    3394 &      6,4 &      4,28 &     2663 &      6,1 &      4,28 &    2192 &      5,1 &      4,28 &     2574 &      5,7 \\

\multicolumn{ 1}{|c|}{} &         25 &          3 &      17,29 &    1218 &      29,0 &      16,95 &     1294 &      31,7 &      16,95 &    1314 &      33,9 &      16,95 &     1167 &      28,5 \\

\multicolumn{ 1}{|c|}{} &         25 &          8 &      7,91 &    6841 &      33,31 &      7,74 &     5315 &      30,7 &      7,74 &    3127 &      22,7 &      7,74 &     5038 &      30,4 \\

\multicolumn{ 1}{|c|}{} &         25 &         10 &       9,81 &   25390 &      102,1 &       9,76 &    11640 &      55,2 &       9,76 &   12560 &      61,8 &       9,76 &    13830 &      66,2 \\

\multicolumn{ 1}{|c|}{} &         28 &          3 &      17,36 &    7432 &      168,4 &       17,10 &     6343 &      155,4 &       17,10 &    7032 &      164,6 &       17,10 &     7531 &      168,8 \\

\multicolumn{ 1}{|c|}{} &         28 &          8 &      10,67 &   34800 &      271,1 &      10,63 &    23080 &      204,8 &      10,63 &   21210 &      191,3 &      10,63 &    20850 &      188,2 \\

\multicolumn{ 1}{|c|}{} &         28 &         10 &      9,141 &  131900 &      832,1 &      9,11 &    92660 &        659,0 &      9,11 &   81600 &      634,5 &      9,11 &    93700 &      705,6 \\

\multicolumn{ 1}{|c|}{} &         30 &          3 &      21,74 &    7220 &      271,6 &      21,47 &     3502 &      164,5 &      21,47 &    4853 &      204,9 &      21,47 &     4873 &      190,1 \\

\multicolumn{ 1}{|c|}{} &         30 &          8 &      14,48 &  273600 &       2896,0 &      14,36 &   183100 &       2022,0 &      14,36 &  204500 &       2250,0 &      14,36 &   258800 &       3026,0 \\

\multicolumn{ 1}{|c|}{} &         30 &         10 &      10,56 &  109700 &      863,2 &       10,50 &    99160 &      837,6 &       10,50 &   94870 &      784,8 &       10,50 &   103700 &      898,1 \\
\cline{2-15}\rowcolor[gray]{0.9}
\multicolumn{ 1}{|c|}{} & \multicolumn{ 2}{c|}{{\bf Average}} & {\bf 12,10} & {\bf 40522} & {\bf 366,3} & {\bf 11,90} & {\bf 29007} & {\bf 279,3} & {\bf 11,90} & {\bf 29298} & {\bf 291,8} & {\bf 11,90} & {\bf 34543} & {\bf 355,4} \\
\hline
\multirow{15}{*}{\begin{sideways}\texttt{TRIMMEAN}\end{sideways}}  &         15 &          3 &      16,44 &      729 &      1,4 &      16,19 &      677 &      1,4 &      16,19 &      733 &      1,4 &      16,19 &      789 &       1,5 \\

\multicolumn{ 1}{|c|}{} &         15 &          5 &      8,58 &      646 &     0,7 &      8,39 &      566 &     0,7 &      8,39 &      550 &      0,8 &      8,39 &      489 &     0,7 \\

\multicolumn{ 1}{|c|}{} &         15 &          8 &      7,88 &      322 &       0,3 &      7,85 &      353 &     0,3 &      7,85 &      349 &      0,3 &      7,85 &      366 &     0,3 \\

\multicolumn{ 1}{|c|}{} &         20 &          3 &      17,65 &    1022 &       6,2 &       17,60 &      871 &      6,5 &       17,60 &      977 &       6,2 &       17,60 &      929 &      6,8 \\

\multicolumn{ 1}{|c|}{} &         20 &          8 &      9,08 &    4512 &      7,9 &      9,00 &     4519 &      8,3 &      9,00 &    4508 &      7,9 &      9,00 &     4476 &      8,4 \\

\multicolumn{ 1}{|c|}{} &         20 &         10 &      5,97 &    1620 &      2,7 &      5,95 &     1402 &      2,5 &      5,95 &    1353 &      2,5 &      5,95 &     1444 &       2,6 \\

\multicolumn{ 1}{|c|}{} &         25 &          3 &      16,87 &    2081 &      25,1 &      16,83 &     1844 &      24,1 &      16,83 &    1879 &      26,8 &      16,83 &     1486 &      22,5 \\

\multicolumn{ 1}{|c|}{} &         25 &          8 &      8,63 &   21460 &      74,2 &      8,40 &    13420 &      54,4 &      8,40 &   14300 &      56,4 &      8,40 &    15390 &      59,5 \\

\multicolumn{ 1}{|c|}{} &         25 &         10 &      12,28 &    1760 &      6,5 &      12,14 &     1106 &      4,8 &      12,14 &    1384 &      5,7 &      12,14 &     1251 &      5,6 \\

\multicolumn{ 1}{|c|}{} &         28 &          3 &      17,49 &    3800 &      69,6 &      17,48 &     4117 &      80,5 &      17,48 &    4005 &      71,0 &      17,48 &     4458 &      80,8 \\

\multicolumn{ 1}{|c|}{} &         28 &          8 &      10,24 &   20290 &      130,9 &      9,97 &    15890 &      111,4 &      9,97 &   17340 &      119,1 &      9,97 &    16910 &      119,1 \\

\multicolumn{ 1}{|c|}{} &         28 &         10 &      9,17 &   14770 &       69,3 &       8,92 &     8775 &      41,6 &       8,92 &    9241 &      44,4 &       8,92 &    10490 &      49,7 \\

\multicolumn{ 1}{|c|}{} &         30 &          3 &      22,08 &    3302 &      124,4 &      22,05 &     4052 &      157,8 &      22,05 &    4356 &        157,0 &      22,05 &     3925 &      147,3 \\

\multicolumn{ 1}{|c|}{} &         30 &          8 &      15,09 &   54920 &      483,1 &      14,99 &    45860 &      450,8 &      14,99 &   50140 &      486,9 &      14,99 &    52040 &      514,7 \\

\multicolumn{ 1}{|c|}{} &         30 &         10 &       9,36 &  107100 &        757,0 &      9,17 &    42140 &      316,1 &      9,17 &   36570 &      275,1 &      9,17 &    36480 &      263,3 \\
\cline{2-15}\rowcolor[gray]{0.9}
\multicolumn{ 1}{|c|}{} & \multicolumn{ 2}{c|}{{\bf Average}} & {\bf 12,50} & {\bf 15889} & {\bf 117,3} & {\bf 12,30} & {\bf 9706} & {\bf 84,1} & {\bf 12,30} & {\bf 9846} & {\bf 84,1} & {\bf 12,30} & {\bf 10062} & {\bf 85,5} \\
\hline
\end{tabular}}
\label{tabla3}
\end{table}

\begin{table}[htbp]
\vspace*{-1.cm}
 \caption{Computational results for Formulation$F_{uv}$ with Valid Inequalities (I) \label{table4}}
 \centering
{
\tiny \setlength{\tabcolsep}{0.2cm}
\vspace*{-0.3cm}
\hspace*{-1cm}

\begin{tabular}{|c|c|c|rrr|rrr|rrr|rrr|}
\hline
  &      \multirow{ 2}{*}{n} & \multirow{ 2}{*}{p}&
\multicolumn{ 3}{c|}{{ P$_{1}$+P$_{2}=$  P$_{12}$ }} &
\multicolumn{ 3}{c|}{{\bf P$_{12}$+V.I.(\ref{in:DV3_3})-(\ref{in:DV3_2})}} &
\multicolumn{ 3}{c|}{{\bf P$_{12}$+V.I.(\ref{in:DV3_3}),(\ref{in:DV3b_1})-(\ref{in:DV3b_2})}} &
\multicolumn{ 3}{c|}{{\bf P$_{12}$+V.I.(\ref{in:DV3_3}),(\ref{in:DV3c_1})-(\ref{in:DV3c_2})}} \\
\cline{4-15}
           &           &          &       RGAP &      Nodes &       time &       RGAP &      Nodes &       time &       RGAP &      Nodes &       time &       RGAP &      Nodes &       time \\
 \hline
\multirow{15}{*}{\begin{sideways}\texttt{CENTRUM}\end{sideways}}  &         15 &          3 &       16,90 &      489 &     0,9 &       16,90 &      429 &     0,9 &       16,90 &      468 &      0,9 &       16,90 &      701 &      1,1 \\

\multicolumn{ 1}{|c|}{} &         15 &          5 &      8,48 &      809 &     0,9 &      8,47 &      788 &     0,9 &      8,47 &      800 &      1,0 &      8,47 &      759 &       1,0 \\

\multicolumn{ 1}{|c|}{} &         15 &          8 &      8,83 &      501 &     0,6 &      8,83 &      495 &     0,6 &      8,83 &      554 &     0,6 &      8,83 &      556 &      0,6 \\

\multicolumn{ 1}{|c|}{} &         20 &          3 &      16,72 &    1017 &      4,8 &       16,70 &      840 &      4,6 &       16,70 &    1028 &      4,7 &       16,70 &      876 &      4,5 \\

\multicolumn{ 1}{|c|}{} &         20 &          8 &      11,58 &    3712 &      6,2 &      11,58 &     3185 &      5,9 &      11,58 &    3778 &      6,6 &      11,58 &     4408 &      7,9 \\

\multicolumn{ 1}{|c|}{} &         20 &         10 &      8,82 &    3019 &      5,3 &      8,82 &     2430 &      4,9 &      8,82 &    2377 &      4,9 &      8,82 &     2044 &      4,7 \\

\multicolumn{ 1}{|c|}{} &         25 &          3 &      14,57 &      938 &      12,3 &      14,56 &      956 &      14,3 &      14,56 &      821 &      12,5 &      14,56 &      983 &      13,6 \\

\multicolumn{ 1}{|c|}{} &         25 &          8 &      11,33 &   26520 &      87,4 &      11,33 &    25490 &      89,2 &      11,33 &   26310 &      95,2 &      11,33 &    30540 &      107,2 \\

\multicolumn{ 1}{|c|}{} &         25 &         10 &      16,16 &    5134 &      17,7 &      16,16 &     3933 &      15,3 &      16,16 &    4357 &      15,9 &      16,16 &     5271 &       19,2 \\

\multicolumn{ 1}{|c|}{} &         28 &          3 &      17,56 &    4172 &      66,8 &      17,55 &     5275 &      75,3 &      17,55 &    4087 &      69,8 &      17,55 &     5190 &       86,6 \\

\multicolumn{ 1}{|c|}{} &         28 &          8 &      15,77 &   29450 &      170,1 &      15,77 &    25760 &      157,3 &      15,77 &   30410 &      185,7 &      15,77 &    25740 &      160,5 \\

\multicolumn{ 1}{|c|}{} &         28 &         10 &      14,53 &   48120 &      218,4 &      14,53 &    49170 &      233,4 &      14,53 &   62630 &        298,0 &      14,53 &    53050 &        269,0 \\

\multicolumn{ 1}{|c|}{} &         30 &          3 &      19,75 &    6536 &      147,3 &      19,73 &     4810 &        124,0 &      19,73 &    6142 &      151,6 &      19,73 &     5994 &      151,6 \\

\multicolumn{ 1}{|c|}{} &         30 &          8 &      14,63 &  103700 &      786,6 &      14,63 &    86370 &      680,3 &      14,63 &  113900 &      904,5 &      14,63 &   108300 &      965,3 \\

\multicolumn{ 1}{|c|}{} &         30 &         10 &      11,77 &   96140 &      593,6 &      11,77 &   128200 &      882,8 &      11,77 &  108700 &      766,9 &      11,77 &    94490 &      675,7 \\
\cline{2-15}\rowcolor[gray]{0.9}
\multicolumn{ 1}{|c|}{} & \multicolumn{ 2}{c|}{{\bf Average}}  & {\bf 13,8} & {\bf 22017} & {\bf 141,3} & {\bf 13,8} & {\bf 22542} & {\bf 152,7} & {\bf 13,80} & {\bf 24424} & {\bf 167,9} & {\bf 13,80} & {\bf 22594} & {\bf 164,6} \\
\hline
\multirow{15}{*}{\begin{sideways}\texttt{K-CENTRUM}\end{sideways}} &         15 &          3 &       17,60 &      588 &       1,8 &      17,58 &      659 &      1,9 &      17,58 &      707 &      2,0 &      17,58 &      783 &      1,9 \\

\multicolumn{ 1}{|c|}{} &         15 &          5 &      7,91 &      498 &      1,0 &      7,85 &      653 &      1,3 &      7,85 &      547 &      1,2 &      7,85 &      635 &      1,2 \\

\multicolumn{ 1}{|c|}{} &         15 &          8 &      4,40 &      185 &     0,4 &      4,33 &      110 &     0,4 &      4,33 &      157 &      0,4 &      4,34 &      162 &     0,4 \\

\multicolumn{ 1}{|c|}{} &         20 &          3 &      18,23 &    1.528 &      9,1 &      18,22 &     1590 &      10,6 &      18,22 &    1287 &      10,0 &      18,22 &     1777 &      10,7 \\

\multicolumn{ 1}{|c|}{} &         20 &          8 &      7,12 &    2.902 &      6,2 &      7,09 &     3290 &      7,8 &      7,09 &    2568 &      6,2 &      7,09 &     2821 &      6,7 \\

\multicolumn{ 1}{|c|}{} &         20 &         10 &      4,04 &    3705 &      5,9 &       4,03 &     2973 &      5,2 &      4,03 &    2769 &       4,8 &      4,03 &     3045 &      5,4 \\

\multicolumn{ 1}{|c|}{} &         25 &          3 &      17,48 &    2693 &      32,3 &      17,47 &     1764 &      28,4 &      17,47 &    1537 &      26,6 &      17,47 &     1639 &      28,0 \\

\multicolumn{ 1}{|c|}{} &         25 &          8 &      7,33 &    3910 &      23,2 &      7,30 &     3239 &      21,0 &      7,30 &    3488 &      24,3 &      7,30 &     3264 &      23,5 \\

\multicolumn{ 1}{|c|}{} &         25 &         10 &      7,34 &   12820 &      52,6 &      7,28 &    12730 &       51,9 &      7,28 &   11770 &      49,6 &      7,28 &    11330 &      48,8 \\

\multicolumn{ 1}{|c|}{} &         28 &          3 &      17,17 &    4971 &      111,3 &      17,16 &     4126 &      103,2 &      17,16 &    4417 &      106,4 &      17,16 &     4146 &      108,7 \\

\multicolumn{ 1}{|c|}{} &         28 &          8 &      9,44 &   21180 &      178,9 &      9,39 &    17390 &      158,4 &      9,39 &   14370 &      140,5 &      9,39 &    14120 &      139,8 \\

\multicolumn{ 1}{|c|}{} &         28 &         10 &      7,63 &   35460 &      228,2 &      7,54 &    39050 &      276,9 &      7,54 &   38760 &      268,8 &      7,54 &    34770 &      247,2 \\

\multicolumn{ 1}{|c|}{} &         30 &          3 &      21,74 &    4588 &      187,6 &      21,72 &     4740 &      188,6 &      21,72 &    4069 &      169,9 &      21,72 &     4437 &      194,9 \\

\multicolumn{ 1}{|c|}{} &         30 &          8 &       12,20 &  171500 &       1830,0 &      12,14 &   149900 &       1737,0 &      12,15 &  218900 &       2710,0 &      12,14 &   114200 &       1343,0 \\

\multicolumn{ 1}{|c|}{} &         30 &         10 &      7,63 &   37330 &        340,0 &      7,54 &    41100 &      392,6 &       7,54 &   34690 &      332,8 &       7,54 &    30780 &      305,4 \\
\cline{2-15}\rowcolor[gray]{0.9}
\multicolumn{ 1}{|c|}{} & \multicolumn{ 2}{c|}{{\bf Average}} & {\bf 11,20} & {\bf 20257} & {\bf 200,6} & {\bf 11,10} & {\bf 18888} & {\bf 199,0} & {\bf 11,10} & {\bf 22669} & {\bf 256,9} & {\bf 11,10} & {\bf 15194} & {\bf 164,4} \\
\hline
\multirow{15}{*}{\begin{sideways}\texttt{BLOCKS}\end{sideways}} &         15 &          3 &      15,97 &    1729 &      3,8 &      15,86 &     2184 &      4,2 &      15,86 &    1481 &      3,7 &      15,86 &     2248 &      4,4 \\

\multicolumn{ 1}{|c|}{} &         15 &          5 &      8,08 &      860 &       2,4 &      7,98 &     1085 &      2,8 &      7,98 &      877 &      2,4 &      7,98 &      877 &      2,7 \\

\multicolumn{ 1}{|c|}{} &         15 &          8 &       6,25 &      715 &      1,1 &      6,23 &      748 &       1,2 &      6,23 &      772 &       1,1 &      6,23 &      726 &      1,1\\

\multicolumn{ 1}{|c|}{} &         20 &          3 &      17,58 &    5777 &       24,7 &      17,41 &     3515 &      19,7 &      17,41 &    4382 &      19,3 &      17,41 &     4007 &      21,0 \\

\multicolumn{ 1}{|c|}{} &         20 &          8 &      7,46 &    8064 &      20,4 &        7,40 &     5603 &      15,8 &        7,40 &    6327 &      16,1 &        7,40 &     5217 &      15,8 \\

\multicolumn{ 1}{|c|}{} &         20 &         10 &      4,00 &    4585 &      9,0 &      3,98 &     4308 &      8,6 &      3,98 &    4091 &      8,5 &      3,98 &     3505 &      7,8 \\

\multicolumn{ 1}{|c|}{} &         25 &          3 &      16,45 &    4087 &      59,1 &      16,28 &     4489 &      67,7 &      16,28 &    4821 &      65,1 &      16,28 &     4434 &      63,9 \\

\multicolumn{ 1}{|c|}{} &         25 &          8 &      8,12 &   12050 &      89,3 &      8,07 &     8427 &      64,9 &      8,07 &    7338 &      54,5 &      8,07 &     7836 &      64,9 \\

\multicolumn{ 1}{|c|}{} &         25 &         10 &      10,01 &   40870 &      177,5 &      9,99 &    49940 &      199,6 &      9,99 &   36030 &      155,2 &      9,99 &    35520 &      158,8 \\

\multicolumn{ 1}{|c|}{} &         28 &          3 &      17,16 &    8569 &      267,1 &      17,11 &    11070 &        312,0 &      17,11 &    5840 &      220,4 &      17,11 &     5934 &      210,7 \\

\multicolumn{ 1}{|c|}{} &         28 &          8 &       9,32 &   54910 &        745,0 &      9,14 &    47500 &        723,0 &      9,14 &   54830 &      749,3 &      9,14 &    43450 &      580,5 \\

\multicolumn{ 1}{|c|}{} &         28 &         10 &      8,67 &  184700 &       1601,0 &      8,61 &   124900 &       1118,0 &      8,61 &  149500 &       1385,0 &      8,61 &   140900 &       1323,0 \\

\multicolumn{ 1}{|c|}{} &         30 &          3 &      22,39 &   12840 &      703,5 &       22,20 &    10580 &      572,6 &       22,20 &   10520 &      542,8 &       22,20 &     9661 &        517,0 \\

\multicolumn{ 1}{|c|}{} &         30 &          8 &      12,95 &  193300 &       5016,0$^{(3)}$ &      12,76 &   261900 &       6298,0$^{(4)}$ &      12,75 &  197700 &       4858,0$^{(3)}$ &      12,74 &   187400 &       4854,0$^{(3)}$\\

\multicolumn{ 1}{|c|}{} &         30 &         10 &      9,45 &  256600 &       3402,0$^{(1)}$ &      9,29 &   248300 &       3476,0$^{(1)}$ &      9,29 &  215800 &       3277,0$^{(3)}$ &      9,30 &   253000 &       3586,0$^{(1)}$ \\
\cline{2-15}\rowcolor[gray]{0.9}
\multicolumn{ 1}{|c|}{} & \multicolumn{ 2}{c|}{{\bf Average}}  & {\bf 11,60} & {\bf 52644} & {\bf 808,1} & {\bf 11,50} & {\bf 52303} & {\bf 858,9} & {\bf 11,50} & {\bf 46687} & {\bf 757,2} & {\bf 11,50} & {\bf 46981} & {\bf 760,8} \\
\hline \rowcolor[gray]{0.8}
   \multicolumn{ 3}{|c|}{{\bf TOTAL}}  & {\bf 11,80} & {\bf 26543} & {\bf 286,5} & {\bf 11,70} & {\bf 23052} & {\bf 274,7} & {\bf 11,70} & {\bf 23029} & {\bf 270,8} & {\bf 11,70} & {\bf 22470} & {\bf 266,2} \\
\hline
\end{tabular}}

\label{tabla4}
\end{table}

\begin{table}[htbp]
\vspace*{-1.cm}
 \caption{Computational results for Formulation$F_{uv}$ with Valid Inequalities (II) \label{table5}}
 \centering
{
\tiny \setlength{\tabcolsep}{0.2cm}
\vspace*{-0.3cm}
\hspace*{-1cm}

\begin{tabular}{|c|c|c|rrr|rrr|rrr|rrr|}

 \hline
&      \multirow{ 2}{*}{n} & \multirow{ 2}{*}{p}&
\multicolumn{ 3}{c|}{{ P$_{12}$+V.I.(\ref{in:DV3_3}),(\ref{in:DV3c_1})-(\ref{in:DV3c_2})}} &
\multicolumn{ 3}{p{3.2cm}|}{{\bf P$_{12}$+V.I.(\ref{in:DV3_3}),(\ref{in:DV3c_1})-(\ref{in:cortepacking1}) }} &
\multicolumn{ 3}{p{3.2cm}|}{{\bf P$_{12}$+V.I.(\ref{in:DV3_3}),(\ref{in:DV3c_1})-(\ref{const1}) }} &
\multicolumn{ 3}{p{3.2cm}|}{{\bf P$_{12}$+V.I.(\ref{in:DV3_3}),(\ref{in:DV3c_1})-(\ref{const3}) }} \\
\cline{4-15}
           &           &          &       RGAP &      Nodes &       time &       RGAP &      Nodes &       time &       RGAP &      Nodes &       time &       RGAP &      Nodes &       time \\
 \hline

\multirow{15}{*}{\begin{sideways}\texttt{MEDIAN}\end{sideways}} &         15 &          3 &      14,21 &      239 &      1,4 &      7,71 &      168 &      2,3 &      7,71 &      197 &      2,5 &      7,63 &      210 &       2,3 \\

\multicolumn{ 1}{|c|}{} &         15 &          5 &        6,50 &      230 &     0,8 &      4,58 &      122 &      1,4 &      4,58 &      117 &      1,4 &      4,60 &      132 &      1,6 \\

\multicolumn{ 1}{|c|}{} &         15 &          8 &      4,33 &      162 &     0,5 &      4,04 &      166 &     0,7 &      4,04 &      150 &     0,8 &      4,04 &      153 &     0,9 \\

\multicolumn{ 1}{|c|}{} &         20 &          3 &      15,66 &      560 &      7,3 &      7,08 &      220 &      11,9 &      7,08 &      345 &      12,3 &      7,16 &      281 &      12,7 \\

\multicolumn{ 1}{|c|}{} &         20 &          8 &      4,98 &      936 &      3,9 &      4,38 &      262 &      4,4 &      4,38 &      600 &      5,4 &      4,36 &      318 &       4,3 \\

\multicolumn{ 1}{|c|}{} &         20 &         10 &      2,33 &      857 &      2,4 &      2,24 &      354 &      3,2 &      2,21 &      328 &      2,9 &      2,21 &      273 &      2,9 \\

\multicolumn{ 1}{|c|}{} &         25 &          3 &      14,79 &      390 &      19,9 &      7,52 &      429 &      44,5 &      7,512 &      319 &      42,9 &      7,52 &      393 &      42,7 \\

\multicolumn{ 1}{|c|}{} &         25 &          8 &      5,39 &     3301 &      26,5 &      4,03 &      521 &      17,8 &      4,02 &      486 &      18,5 &      4,02 &      405 &      16,3 \\

\multicolumn{ 1}{|c|}{} &         25 &         10 &       5,97 &     2734 &       16,8 &      5,48 &    1177 &      15,9 &      5,48 &     1265 &      16,4 &      5,48 &    1175 &      17,3 \\

\multicolumn{ 1}{|c|}{} &         28 &          3 &      15,72 &     2055 &      71,9 &      6,75 &      320 &      88,1 &      6,75 &      414 &      82,8 &      6,76 &      322 &       85,7 \\

\multicolumn{ 1}{|c|}{} &         28 &          8 &      7,53 &     6398 &      70,4 &      5,15 &    1627 &      47,6 &      5,13 &     1294 &      42,9 &      5,15 &    1781 &      49,3 \\

\multicolumn{ 1}{|c|}{} &         28 &         10 &      6,48 &    17930 &      125,4 &      4,73 &    4394 &      57,9 &      4,73 &     4914 &      62,8 &      4,70 &    3674 &      50,2 \\

\multicolumn{ 1}{|c|}{} &         30 &          3 &       20,60 &     2770 &      142,3 &       13,90 &      868 &      143,6 &      13,67 &      846 &      160,5 &      13,71 &      889 &      157,5 \\

\multicolumn{ 1}{|c|}{} &         30 &          8 &      10,35 &    25060 &      332,7 &      8,19 &    4083 &      111,7 &      8,19 &     3348 &      96,0 &      8,24 &    4593 &      127,6 \\

\multicolumn{ 1}{|c|}{} &         30 &         10 &      6,75 &    18130 &      179,8 &      5,49 &    4232 &      77,7 &      5,49 &     3426 &      65,9 &      5,49 &    3.64 &      76,5 \\
\cline{2-15}\rowcolor[gray]{0.9}
\multicolumn{ 1}{|c|}{} & \multicolumn{ 2}{c|}{{\bf Average}} &  {\bf 9,40} & {\bf 5450} & {\bf 66,8} &  {\bf 6,10} & {\bf 1263} & {\bf 41,9} &  {\bf 6,10} & {\bf 1203} & {\bf 40,9} &  {\bf 6,10} & {\bf 1231} & {\bf 43,2} \\
\hline
\multirow{15}{*}{\begin{sideways}\texttt{ANTI-TRIMMEAN}\end{sideways}} &         15 &          3 &      17,18 &      556 &      2,0 &      15,28 &      467 &      2,9 &      15,28 &      526 &      3,3 &      15,28 &      510 &      2,9 \\

\multicolumn{ 1}{|c|}{} &         15 &          5 &      8,36 &      705 &      1,5 &      8,02 &      562&        1,9 &      8,02 &      432 &      1,7 &      8,02 &      441 &      1,8 \\

\multicolumn{ 1}{|c|}{} &         15 &          8 &       6,82 &      736 &     1,0&      6,71 &      669 &      1,3 &      6,71 &      686 &       1,2 &      6,71 &      733 &      1,3 \\

\multicolumn{ 1}{|c|}{} &         20 &          3 &         17,00 &     1370 &      10,9 &      14,55 &      737 &      15,6 &      14,55 &     1027 &      15,9 &      14,55 &    1008 &      16,6 \\

\multicolumn{ 1}{|c|}{} &         20 &          8 &      7,54 &     2709 &      7,7 &      7,31 &    1876 &      7,7 &      7,31 &     1684 &      7,1 &      7,31 &    2041 &      7,9 \\

\multicolumn{ 1}{|c|}{} &         20 &         10 &      4,28 &     2574 &      5,7 &       4,21 &    2168 &      6,2 &       4,21 &     2367 &      6,4 &       4,21 &    2720 &      7,1 \\

\multicolumn{ 1}{|c|}{} &         25 &          3 &      16,95 &     1167 &      28,5 &      15,41 &    1016 &      47,5 &      15,41 &     1055 &      50,6 &      15,41 &      800 &      50,9 \\

\multicolumn{ 1}{|c|}{} &         25 &          8 &      7,74 &     5038 &      30,4 &      7,03 &    1903 &       22,9 &      7,03 &     2066 &       23,7 &      7,03 &    2664 &      26,7 \\

\multicolumn{ 1}{|c|}{} &         25 &         10 &       9,76 &    13830 &      66,2 &      9,69 &   13360 &      71,2 &      9,69 &    10670 &      57,1 &      9,69 &   13710 &      70,3 \\

\multicolumn{ 1}{|c|}{} &         28 &          3 &       17,10 &     7531 &      168,8 &      14,25 &    2024 &      112,8 &      14,25 &     1943 &        114,0 &      14,25 &    2380 &      119,8 \\

\multicolumn{ 1}{|c|}{} &         28 &          8 &      10,63 &    20850 &      188,2 &      10,34 &   11240 &      133,7 &      10,34 &    12110 &      144,9 &      10,34 &   10960 &        128,0 \\

\multicolumn{ 1}{|c|}{} &         28 &         10 &      9,11 &    93700 &      705,6 &       8,89 &   66920 &      519,5 &       8,89 &    61180 &      522,6 &       8,89 &   69100 &      527,7 \\

\multicolumn{ 1}{|c|}{} &         30 &          3 &      21,47 &     4873 &      190,1 &       18,10 &    4082 &      295,5 &       18,10 &     2724 &        255,0 &       18,10 &    2702 &      239,1 \\

\multicolumn{ 1}{|c|}{} &         30 &          8 &      14,36 &   258800 &       3026,0 &      13,96 &  135200 &       1701,0 &      13,96 &   124200 &       1576,0 &      13,96 &  123600 &       1526,0 \\

\multicolumn{ 1}{|c|}{} &         30 &         10 &       10,50 &   103700 &      898,1 &      10,26 &   55810 &      497,4 &      10,26 &    61050 &      566,9 &      10,33 &   65660 &      606,6 \\
\cline{2-15}\rowcolor[gray]{0.9}
\multicolumn{ 1}{|c|}{} & \multicolumn{ 2}{c|}{{\bf Average}} & {\bf 11,90} & {\bf 34543} & {\bf 355,4} & {\bf 10,90} & {\bf 19869} & {\bf 229,1} & {\bf 10,90} & {\bf 18915} & {\bf 223,1} & {\bf 10,90} & {\bf 19935} & {\bf 222,2} \\
\hline
\multirow{15}{*}{\begin{sideways}\texttt{TRIMMEAN}\end{sideways}} &         15 &          3 &      16,19 &      789 &       1,5 &      14,68 &      647 &      2,1 &      14,56 &      546 &      2,3 &      14,56 &      437&      2,2 \\

\multicolumn{ 1}{|c|}{} &         15 &          5 &      8,39 &      489 &     0,7 &      7,71 &      525 &      1,1 &      7,71 &      487 &      1,0 &      7,71 &      513 &      1,1 \\

\multicolumn{ 1}{|c|}{} &         15 &          8 &      7,85 &      366 &     0,3 &      7,79 &      342 &     0,5 &      7,79 &      417 &      0,6 &      7,79 &      384 &     0,5 \\

\multicolumn{ 1}{|c|}{} &         20 &          3 &       17,60 &      929 &      6,8 &      15,17 &      978 &      11,8 &      15,02 &      648 &      11,7 &      15,02 &      861 &      12,8 \\

\multicolumn{ 1}{|c|}{} &         20 &          8 &      9,01 &     4476 &      8,4 &       8,68 &    4139 &      8,5 &      8,59 &     3603 &      8,5 &      8,59 &    3434 &       8,5 \\

\multicolumn{ 1}{|c|}{} &         20 &         10 &      5,95 &     1444 &       2,6&      5,83 &    1302 &      3,6 &      5,80 &     1251 &      3,8 &      5,80 &    1464 &      4,2 \\

\multicolumn{ 1}{|c|}{} &         25 &          3 &      16,83 &     1486 &      22,5 &      15,84 &    1671 &       39,2 &      15,74 &     1077 &      38,9 &      15,74 &    1148 &      39,4 \\

\multicolumn{ 1}{|c|}{} &         25 &          8 &      8,40 &    15390 &      59,5 &       8,24 &   12770 &      57,5 &      8,19 &    11340 &      53,7 &      8,19 &   11400 &      54,1 \\

\multicolumn{ 1}{|c|}{} &         25 &         10 &      12,14 &     1251 &      5,6 &      12,07 &    1382 &       9,5 &      12,07 &     1249 &      8,9 &      12,07 &    1155 &      9,0 \\

\multicolumn{ 1}{|c|}{} &         28 &          3 &      17,48 &     4458 &      80,8 &      15,77 &    2371 &      88,3 &       15,60 &     1727 &      83,5 &       15,60 &    1721 &      82,9 \\

\multicolumn{ 1}{|c|}{} &         28 &          8 &      9,97 &    16910 &      119,1 &      9,63 &   11070 &      93,3 &      9,39 &     9072 &      81,6 &      9,39 &    7802 &      74,6 \\

\multicolumn{ 1}{|c|}{} &         28 &         10 &       8,92 &    10490 &      49,7 &      8,76 &    9234 &      55,4 &      8,76 &    10090 &      57,1 &      8,76 &   10920 &      63,8 \\

\multicolumn{ 1}{|c|}{} &         30 &          3 &      22,05 &     3925 &      147,3 &       19,80 &    3477 &        200,0 &      19,69 &     2634 &      189,6 &      19,69 &    2233 &        177,0 \\

\multicolumn{ 1}{|c|}{} &         30 &          8 &      14,99 &    52040 &      514,7 &      14,54 &   33800 &      348,4 &      14,52 &    29400 &      326,1 &      14,52 &   29980 &      339,6 \\

\multicolumn{ 1}{|c|}{} &         30 &         10 &      9,17 &    36480 &      263,3 &       8,63 &   26750 &      205,3 &       8,63 &    26950 &      229,9 &       8,63 &   27860 &      228,9 \\
\cline{2-15}\rowcolor[gray]{0.9}
\multicolumn{ 1}{|c|}{} & \multicolumn{ 2}{c|}{{\bf Average}} & {\bf 12,30} & {\bf 10062} & {\bf 85,5} & {\bf 11,50} & {\bf 7364} & {\bf 75,0} & {\bf 11,50} & {\bf 6699} & {\bf 73,1} & {\bf 11,50} & {\bf 6754} & {\bf 73,3} \\
\hline
\end{tabular}
}

\label{tabla5}
\end{table}

\begin{table}[htbp]
\vspace*{-1.cm}
 \caption{Computational results for Formulation$F_{uv}$ with Valid Inequalities (II){\it bis} \label{table6}}
 \centering
{
\tiny \setlength{\tabcolsep}{0.2cm}
\vspace*{-0.3cm}
\hspace*{-1cm}

\begin{tabular}{|c|c|c|rrr|rrr|rrr|rrr|}

 \hline
&      \multirow{ 2}{*}{n} & \multirow{ 2}{*}{p}&
\multicolumn{ 3}{c|}{{ P$_{12}$+V.I.(\ref{in:DV3_3}),(\ref{in:DV3c_1})-(\ref{in:DV3c_2})}} &
\multicolumn{ 3}{p{3.2cm}|}{{\bf P$_{12}$+V.I.(\ref{in:DV3_3}),(\ref{in:DV3c_1})-(\ref{in:cortepacking1}) }} &
\multicolumn{ 3}{p{3.2cm}|}{{\bf P$_{12}$+V.I.(\ref{in:DV3_3}),(\ref{in:DV3c_1})-(\ref{const1}) }} &
\multicolumn{ 3}{p{3.2cm}|}{{\bf P$_{12}$+V.I.(\ref{in:DV3_3}),(\ref{in:DV3c_1})-(\ref{const3}) }} \\
\cline{4-15}
           &           &          &       RGAP &      Nodes &       time &       RGAP &      Nodes &       time &       RGAP &      Nodes &       time &       RGAP &      Nodes &       time \\
 \hline

\multirow{15}{*}{\begin{sideways}\texttt{CENTRUM}\end{sideways}} &         15 &          3 &       16,90 &      701&      1,1 &      16,32 &      574 &      1,6 &      16,32 &      576 &      1,6 &      16,32 &      513 &      1,6 \\

\multicolumn{ 1}{|c|}{} &         15 &          5 &      8,47 &      759 &       1,0 &      8,45 &      769 &      1,3 &      8,45 &      850 &      1,4 &      8,45 &      850 &      1,5 \\

\multicolumn{ 1}{|c|}{} &         15 &          8 &      8,83 &      556 &      0,6 &      8,81 &      577 &     0,8 &      8,78 &      582 &     0,9 &      8,78 &      437 &     0,9 \\

\multicolumn{ 1}{|c|}{} &         20 &          3 &       16,70 &      876 &      4,5 &      15,94 &      579 &      7,1 &      15,94 &      649 &      7,0 &      15,94 &      705 &      6,9 \\

\multicolumn{ 1}{|c|}{} &         20 &          8 &      11,58 &     4408 &      7,9 &      11,57 &    4061 &      7,9 &      11,56 &     3902 &      8,3 &      11,56 &    3725 &      7,8 \\

\multicolumn{ 1}{|c|}{} &         20 &         10 &      8,82 &     2044 &      4,7 &      8,82 &    2044 &      5,4 &      8,82 &     2739 &      6,2 &      8,82 &    2389 &      5,6 \\

\multicolumn{ 1}{|c|}{} &         25 &          3 &      14,56 &      983 &      13,6 &      14,27 &      689 &      19,0 &      14,27 &      669 &      19,8 &      14,27 &      750 &       19,9 \\

\multicolumn{ 1}{|c|}{} &         25 &          8 &      11,33 &    30540 &      107,2 &      11,29 &   29000 &      103,6 &      11,29 &    25290 &      96,5 &      11,29 &   28560 &      104,6 \\

\multicolumn{ 1}{|c|}{} &         25 &         10 &      16,16 &     5271 &       19,2 &      16,15 &    4339 &      18,4 &      16,15 &     4335 &      18,6 &      16,15 &    5226 &      21,4 \\

\multicolumn{ 1}{|c|}{} &         28 &          3 &      17,55 &     5190 &       86,6 &      16,87 &    2660 &      74,6 &      16,87 &     2788 &      79,9 &      16,87 &    2997 &      83,0 \\

\multicolumn{ 1}{|c|}{} &         28 &          8 &      15,77 &    25740 &      160,5 &      15,75 &   26910 &      175,6 &      15,75 &    25320 &      179,3 &      15,75 &   26240 &      173,8 \\

\multicolumn{ 1}{|c|}{} &         28 &         10 &      14,53 &    53050 &        269,0 &      14,52 &   57390 &      296,1 &      14,52 &    58690 &      314,7 &      14,52 &   50760 &      265,1 \\

\multicolumn{ 1}{|c|}{} &         30 &          3 &      19,73 &     5994 &      151,6 &      18,36 &    4419 &      183,1 &      18,36 &     3236 &        148,0 &      18,36 &    4846 &      190,2 \\

\multicolumn{ 1}{|c|}{} &         30 &          8 &      14,63 &   108300 &      965,3 &      14,59 &  119500 &       1068,0 &      14,59 &   110300 &       1094,0 &      14,59 &  133400 &       1228,0 \\

\multicolumn{ 1}{|c|}{} &         30 &         10 &      11,77 &    94490 &      675,7 &      11,77 &   94490 &      678,4 &      11,77 &    91310 &      700,3 &      11,77 &  107000 &      806,9 \\
\cline{2-15}\rowcolor[gray]{0.9}
\multicolumn{ 1}{|c|}{} & \multicolumn{ 2}{c|}{{\bf Average}} & {\bf 13,80} & {\bf 22594} & {\bf 164,6} & {\bf 13,60} & {\bf 23200} & {\bf 176,1} & {\bf 13,60} & {\bf 22082} & {\bf 178,4} & {\bf 13,60} & {\bf 24560} & {\bf 194,5} \\
\hline
\multirow{15}{*}{\begin{sideways}\texttt{K-CENTRUM}\end{sideways}} &         15 &          3 &      17,58 &      783 &      1,9 &      14,72 &      678 &      3,3 &      14,49 &      502 &      3,2 &      14,49 &      594 &      3,4 \\

\multicolumn{ 1}{|c|}{} &         15 &          5 &      7,85 &      6345 &      1,2 &      6,80 &      437 &      1,8 &      6,58 &      387 &      1,8 &      6,58 &      358 &      1,9 \\

\multicolumn{ 1}{|c|}{} &         15 &          8 &      4,34 &      162 &     0,4 &      4,04 &      166 &     0,7 &      4,04 &      150 &     0,8 &      4,04 &      153 &     0,8 \\

\multicolumn{ 1}{|c|}{} &         20 &          3 &      18,22 &     1777 &      10,7 &      15,52 &      959 &      13,3 &      15,25 &      847 &      14,1 &      15,25 &      904 &      13,8 \\

\multicolumn{ 1}{|c|}{} &         20 &          8 &      7,09 &     2821 &      6,7 &      6,77 &    1272 &      6,1 &      6,62 &     1260 &      6,7 &      6,62 &    1001 &      6,4 \\

\multicolumn{ 1}{|c|}{} &         20 &         10 &      4,03 &     3045 &      5,4 &       3,84 &    2437 &      6,0 &      3,74 &     1879 &      5,6 &      3,74 &    1950 &      5,5 \\

\multicolumn{ 1}{|c|}{} &         25 &          3 &      17,47 &     1639 &      28,0 &      16,17 &    1404 &      46,7 &      16,11 &     1035 &      45,3 &      16,11 &    1377 &      50,7 \\

\multicolumn{ 1}{|c|}{} &         25 &          8 &      7,29 &     3264 &      23,5 &      6,46 &    1895 &      23,2 &      6,30 &     1777 &      23,8 &      6,30 &    1571 &      23,1 \\

\multicolumn{ 1}{|c|}{} &         25 &         10 &      7,28 &    11330 &      48,8 &      7,09 &    6548 &      38,5 &      6,99 &     5529 &      36,0 &      6,99 &    5793 &      37,2 \\

\multicolumn{ 1}{|c|}{} &         28 &          3 &      17,16 &     4146 &      108,7 &      14,52 &    1557 &       83,6 &      14,21 &     1076 &      89,4 &      14,21 &      904 &      82,9 \\

\multicolumn{ 1}{|c|}{} &         28 &          8 &      9,39 &    14120 &      139,8 &      8,345 &    5935 &      89,7 &      8,14 &     5073 &      82,3 &      8,14 &    5264 &      84,9 \\

\multicolumn{ 1}{|c|}{} &         28 &         10 &      7,54 &    34770 &      247,2 &      6,66 &   26570 &      206,3 &      6,51 &    11410 &      113,9 &      6,51 &    9299 &      99,6 \\

\multicolumn{ 1}{|c|}{} &         30 &          3 &      21,72 &     4437 &      194,9 &      17,34 &    3598 &      269,3 &      17,02 &     2889 &      237,2 &      17,02 &    3989 &      279,3 \\

\multicolumn{ 1}{|c|}{} &         30 &          8 &      12,14 &   114200 &       1343,0 &      10,75 &   17250 &        291,0 &      10,66 &    97400 &       1434,0 &      10,77 &   16090 &      303,5 \\

\multicolumn{ 1}{|c|}{} &         30 &         10 &       7,54 &    30780 &      305,4 &      6,60 &   13560 &      170,7 &      6,55 &     8198 &      120,1 &      6,59 &    8883 &      133,7 \\
\cline{2-15}\rowcolor[gray]{0.9}
\multicolumn{ 1}{|c|}{} & \multicolumn{ 2}{c|}{{\bf Average}} & {\bf 11,10} & {\bf 15194} & {\bf 164,4} &  {\bf 9,70} & {\bf 5618} & {\bf 83,3} &  {\bf 9,50} & {\bf 9294} & {\bf 147,6} &  {\bf 9,60} & {\bf 3875} & {\bf 75,1} \\
\hline
\multirow{15}{*}{\begin{sideways}\texttt{BLOCKS}\end{sideways}} &         15 &          3 &      15,86 &     2248 &      4,4 &      13,58 &    2925 &      6,9 &      13,55 &     2611 &      6,7 &      13,55 &    2388 &       6,7 \\

\multicolumn{ 1}{|c|}{} &         15 &          5 &      7,98 &      877 &      2,7 &      7,35 &      785 &      3,5 &      7,35 &      898 &      3,7 &      7,35 &      762 &      3,4 \\

\multicolumn{ 1}{|c|}{} &         15 &          8 &      6,23 &      726 &      1,1 &      6,14 &      678 &      1,4 &      6,11 &      804 &      1,6 &      6,11 &      698 &      1,5 \\

\multicolumn{ 1}{|c|}{} &         20 &          3 &      17,41 &     4007 &      21,0 &      14,06 &    4384 &      31,1 &      14,06 &     4642 &      31,9 &      14,06 &    4712 &      32,7 \\

\multicolumn{ 1}{|c|}{} &         20 &          8 &        7,40 &     5217 &      15,8 &      7,12 &    4646 &      15,9 &       7,11 &     5109 &      17,8 &       7,11 &    4380 &      15,7 \\

\multicolumn{ 1}{|c|}{} &         20 &         10 &      3,98 &     3505 &      7,8 &      3,91 &    2882 &      8,5 &      3,89 &     3016 &      8,8 &      3,89 &    3319 &      9,5 \\

\multicolumn{ 1}{|c|}{} &         25 &          3 &      16,28 &     4434 &      63,9 &      14,65 &    5324 &      98,9 &      14,65 &     4923 &      91,3 &      14,65 &    6696 &        102,0 \\

\multicolumn{ 1}{|c|}{} &         25 &          8 &      8,07 &     7836 &      64,9 &      7,39 &    4695 &      51,8 &      7,39 &     5066 &      54,5 &      7,39 &    5400 &      57,9 \\

\multicolumn{ 1}{|c|}{} &         25 &         10 &      9,99 &    35520 &      158,8 &      9,96 &   34470 &      151,7 &      9,96 &    49260 &      222,5 &      9,96 &   37030 &      169,3 \\

\multicolumn{ 1}{|c|}{} &         28 &          3 &      17,11 &     5934 &      210,7 &      14,05 &    4436 &        229,0 &      14,05 &     5004 &      290,9 &      14,05 &    3642 &        206,0 \\

\multicolumn{ 1}{|c|}{} &         28 &          8 &      9,14 &    43450 &      580,5 &      8,61 &   54300 &       1787,0 &      8,61 &    35590 &      559,7 &      8,61 &   32250 &      503,8 \\

\multicolumn{ 1}{|c|}{} &         28 &         10 &      8,61 &   140900 &       1323,0 &      8,20 &   94410 &       1013,0 &      8,20 &   105800 &       1150,0 &      8,20 &  125300 &       1328,0 \\

\multicolumn{ 1}{|c|}{} &         30 &          3 &       22,20 &     9661 &        517,0 &      18,45 &    7648 &      623,9 &      18,45 &     7766 &      637,6 &      18,45 &    6555 &      539,6 \\

\multicolumn{ 1}{|c|}{} &         30 &          8 &      12,74 &   187400 &       4854,0$^{(3)}$ &      12,26 &  177000 &       4807,0$^{(3)}$ &      12,17 &   176000 &       4400,0$^{(2)}$ &      12,19 &  184600 &       4494,0$^{(2)}$\\

\multicolumn{ 1}{|c|}{} &         30 &         10 &      9,30 &   253000 &       3586,0$^{(1)}$ &      9,06 &  191400 &       2931,0$^{(1)}$ &      9,04 &   208800 &       3333,0$^{(1)}$ &      9,04 &  194700 &       2991,0$^{(1)}$ \\
\cline{2-15}\rowcolor[gray]{0.9}
\multicolumn{ 1}{|c|}{} & \multicolumn{ 2}{c|}{{\bf Average}} & {\bf 11,50} & {\bf 46981} & {\bf 760,8} & {\bf 10,30} & {\bf 39332} & {\bf 784,0} & {\bf 10,30} & {\bf 41019} & {\bf 720,7} & {\bf 10,30} & {\bf 40829} & {\bf 697,4} \\
\hline \rowcolor[gray]{0.8}
     \multicolumn{ 3}{|c|}{{\bf TOTAL}} & {\bf 11,70} & {\bf 22470} & {\bf 266,2} & {\bf 10,40} & {\bf 16108} & {\bf 231,6} & {\bf 10,30} & {\bf 16536} & {\bf 230,7} & {\bf 10,30} & {\bf 16197} & {\bf 217,6} \\
\hline
\end{tabular}
}

\label{tabla6}
\end{table}


\end{document}